\documentclass[12pt]{amsart}

\usepackage{amscd}
\usepackage{amsthm}
\usepackage{amsmath}
\usepackage{latexsym}
\usepackage{graphicx}
\usepackage{delarray}
\usepackage{float}
\usepackage{mathrsfs}
\usepackage[psamsfonts]{amssymb}
\usepackage{vmargin}
\usepackage{amsmath}
\usepackage{epsfig,verbatim}
\usepackage{amsmath,amsfonts,amsthm,amsopn,mathrsfs}
\usepackage{amsfonts,amsthm,amsopn}
\usepackage{float}

\begin{document}


\newtheorem{theorem}{Theorem}[section]
\newtheorem{definition}[theorem]{Definition}
\newtheorem{lemma}[theorem]{Lemma}
\newtheorem{proposition}[theorem]{Proposition}
\newtheorem{corollary}[theorem]{Corollary}
\newtheorem{example}[theorem]{Example}
\newtheorem{remark}[theorem]{Remark}



\def\noo{\noindent }
\newcommand{\rmf}{\rm f}
\newcommand{\rmg}{\rm g}
\newcommand{\rmh}{\rm h}
\newcommand{\Bff}{{\bf f}}
\newcommand{\Bfg}{{\bf g}}
\newcommand{\Bfh}{{\bf h}}
\newcommand{\rmaj}{\rm a(j)}
\newcommand{\rmakr}{\rm a(k,r)}
\newcommand{\rmals}{\rm a(l,s)}
\newcommand{\rmala}{\rm a(\lambda)}
\newcommand{\rmbj}{\rm b(j)}
\newcommand{\rmbkr}{\rm b(k,r)}
\newcommand{\rmbls}{\rm b(l,s)}
\newcommand{\rmbla}{\rm b(\lambda)}
\newcommand{\rmfj}{\rm f(j)}
\newcommand{\rmk}{\rm k}
\newcommand{\rml}{\rm l}
\newcommand{\rmr}{\rm r}
\newcommand{\rms}{\rm s}
\newcommand{\rmM}{\rm M}
\newcommand{\rmN}{\rm N}
\newcommand{\rmZZN}{\mathbb{Z}_{\rm N}}
\newcommand{\rmdZZN}{\widehat{\mathbb{Z}}_{\rm N}}
\newcommand{\rmZZNN}{\mathbb{Z}_{\rm N}\times\rmdZZN}
\newcommand{\rmG}{\rm G}
\newcommand{\rmGG}{\rm G\times\widehat G}
\newcommand{\rmCCM}{\mathbb{C}^{\rm M}}
\newcommand{\rmCCN}{\mathbb{C}^{\rm N}}
\newcommand{\rmCCG}{\mathbb{C}^{|\rmG|}}
\newcommand{\rmCCGG}{\mathbb{C}^{|\rmGG|}}
\newcommand{\rmpikr}{\pi({\rm k},{\rm r})}
\newcommand{\rmpils}{\pi({\rm l},{\rm s})}
\newcommand{\rmA}{\rm A}
\newcommand{\rmB}{\rm B}
\newcommand{\rmIN}{\rm I_N}
\newcommand{\rmC}{\rm C}
\newcommand{\rmCgLa}{{\rm C}_{\Bfg,\Lambda}}
\newcommand{\rmDgLa}{{\rm D}_{\Bfg,\Lambda}}
\newcommand{\rmGgLa}{{\rm G}_{\Bfg,\Lambda}}
\newcommand{\rmR}{\rm R}
\newcommand{\rmS}{\rm S}
\newcommand{\rmSgLa}{{\rm S}_{\Bfg,\Lambda}}
\newcommand{\rmVgf}{\rm V_\Bfg\Bff}
\newcommand{\rmVgfkr}{\rm V_\Bfg\Bff(\rmk,\rmr)}
\newcommand{\rmVgfla}{\rm V_\Bfg\Bff(\lambda)}
\newcommand{\rmVgg}{\rm V_\Bfg\Bfg}
\newcommand{\rmVggkr}{\rm V_\Bfg\Bfg(\rmk,\rmr)}
\newcommand{\rmVggla}{\rm V_\Bfg\Bfg(\lambda)}
\newcommand{\MMNC}{{\mathcal M}_{\rmM\times\rmN}(\CC)}
\newcommand{\MMC}{{\mathcal M}_{\rmM\times\rmM}(\CC)}
\newcommand{\MNC}{{\mathcal M}_{\rmN\times\rmN}(\CC)}
\newcommand{\MMG}{{\mathcal M}_{|\rm G\times\rm G|}(\CC)}
\newcommand{\MMGG}{{\mathcal M}_{|\rm G\times\widehat{\rm G}|}(\CC)}
\newcommand{\ulrmk}{{\rm k}}
\newcommand{\ulrml}{{\rm l}}
\newcommand{\ulrmr}{{\rm r}}
\newcommand{\ulrms}{{\rm s}}
\newcommand{\ulrmaj}{\rm a({j})}
\newcommand{\ulrmakr}{\rm a({k},{r})}
\newcommand{\ulrmals}{\rm a({l},{s})}
\newcommand{\ulrmbj}{\rm b({j})}
\newcommand{\ulrmbkr}{\rm b({k},{r})}
\newcommand{\ulrmbls}{\rm b({l},{s})}
\newcommand{\ulrmpikr}{\pi({\rm k},{\rm r})}
\newcommand{\ulrmpils}{\pi({\rm l},{\rm s})}

\newcommand{\C}{\mathbb{C}}
\newcommand{\CC}{\mathbb{C}}
\newcommand{\CCM}{\mathbb{C}^M}
\newcommand{\CCN}{\mathbb{C}^N}

\newcommand{\NN}{\mathbb{N}}
\newcommand{\TT}{\mathbb{T}}
\newcommand{\ZZ}{\mathbb{Z}}
\newcommand{\phihat}{\widehat{\phi}}
\newcommand{\lahat}{\widehat{\la}}
\newcommand{\fhat}{\widehat{f}}
\newcommand{\fkhat}{\widehat{f_k}}
\newcommand{\La}{\Lambda}
\newcommand{\Om}{\Omega}
\newcommand{\Lap}{\La^\perp}
\newcommand{\la}{\lambda}
\newcommand{\dLa}{{\La^\circ}}
\newcommand{\dla}{{\la^\circ}}
\newcommand{\al}{\alpha}
\newcommand{\be}{\beta}
\newcommand{\tal}{\tilde{\alpha}}
\newcommand{\FF}{\mathcal{F}}
\newcommand{\M}{\mathcal{M}}
\newcommand{\Fro}{\operatorname{Fro}}
\newcommand{\om}{\omega}
\newcommand{\ol}{\overline}
\newcommand{\A}{\mathcal A}
\newcommand{\B}{\mathcal B}
\newcommand{\Bh}{\B(\cH)}
\newcommand{\cC}{\mathcal C}
\newcommand{\cH}{\mathcal H}
\newcommand{\CS}{C^*}
\newcommand{\LtRd}{L^2(\mathbb{R}^d)}
\newcommand{\abs}[1]{\lvert#1\rvert}
\newcommand{\tr}[1]{\operatorname{tr}(#1)}
\newcommand{\absbig}[1]{\big\lvert#1\big\rvert}
\newcommand{\norm}[1]{\lVert#1\rVert}
\newcommand{\norminf}[1]{\norm{#1}_\infty}
\newcommand{\scp}[1]{\langle#1\rangle}
\newcommand{\spn}[1]{\operatorname{span}(#1)}
\newcommand{\scpbig}[1]{\big\langle#1\big\rangle}
\newcommand{\scpBig}[1]{\Big\langle#1\Big\rangle}
\newcommand{\set}[2]{\big\{ \, #1 \, \big| \, #2 \, \big\}}
\newcommand{\inner}[2]{\langle #1,#2\rangle}
\newcommand{\spec}[1]{\operatorname{Sp}(#1)}


\title{{G}abor {a}nalysis {o}ver finite {A}belian groups}

\author{Hans G. Feichtinger}
\address{Fakult\"at f\"ur Mathematik\\
Universit\"at Wien\\
Nordbergstrasse 15\\ 1090 Wien\\ Austria}
\email{hans.feichtinger@univie.ac.at}

\author{Werner Kozek} 
\address{Supremum Solutions\\ Vienna\\ Austria} 
\email{werner.kozek@supremum.at}

\author{Franz Luef$^1$} \footnotetext[1]{F.L. was supported by the Marie-Curie Excellence Grant MEXT-CT-2004-517154.}
\address{Fakult\"at f\"ur Mathematik\\
Universit\"at Wien\\
Nordbergstrasse 15\\ 1090 Wien\\ Austria} \email{franz.luef@univie.ac.at}

\vspace{-16mm}
\begin{abstract}  Gabor frames for signals over
{\it finite Abelian groups}, generated by an arbitrary lattice
within the finite time-frequency plane, are the central topic of this paper. Our generic approach  covers
both multi-dimensional signals as well as non-separable lattices, and in fact the multi-window case as well.
Our generic approach includes most of the fundamental facts about Gabor
expansions of finite signals for the case of product lattices, as they have
been given by Qiu, Wexler-Raz or Tolimieri-Orr, Bastiaans and Van-Leest and others.
In our presentation the  spreading representation of linear
operators between finite-dimensional Hilbert space as well as  a symplectic version of Poisson's summation formula over the finite
time-frequency plane are essential ingredients. They bring us to the so-called Fundamental Identity of Gabor Analysis.
In addition, we highlight projective representations of the time-frequency
plane and its subgroups and explain the natural connection to twisted group
algebras. In the finite-dimensional setting discussed in this paper these twisted group algebras are
just matrix algebras and their structure provides the algebraic framework for
the study of the deeper properties of finite-dimensional Gabor frames, independent of the structure theory 
theorem for finite Abelian groups. 
\end{abstract}

\maketitle \pagestyle{myheadings} \markboth{H.G. Feichtinger, W. Kozek and F.
Luef}{Gabor analysis over finite abelian groups}

\section{Introduction}
In the last two decades a new branch of time-frequency analysis, called {\it
Gabor analysis}, has found many applications in pure and applied mathematics.
The connection between time-variant systems in communication theory and Gabor
analysis has turned out to be of great importance for the application of
Gabor analysis to real-world problems such as the transmission of signals
between cellular phones. In modern digital communication there is an ongoing
trend towards FFT-based multicarrier modulation: popular wireline systems for
the internet access such as ADSL are based on discrete multitone-modulation
(DMT) and important wireless systems such as WLAN, UMTS, WIMAX make use of
OFDM (orthogonal frequency division multiplex)-type modulation. In
mathematical terminology both DMT and OFDM are essentially Weyl-Heisenberg
group structured Riesz bases. Thanks to the Wexler-Raz relation all the mathematical machinery
developed for Gabor frame theory can be exploited for the advanced design of
so-called pulse-shaped multicarrier systems, \cite{st06}. All these
applications lead in a natural manner to the discussion of Gabor frames for
finite-dimensional Hilbert spaces. Due to the large potential of the material
to the above mentioned applications of Gabor frames in signal analysis we
want to address this note not only to mathematicians but also to engineers
who are interested in the deeper mathematical background which is likely to
provide a sound basis for further applications. For the same reason we have
tried to
present the results in a self-contained, although very general form. Finally
we want to mention that many researchers have contributed to the construction
of frames for finite-dimensional Hilbert spaces, e.g.
\cite{befi03,boel03,hest03}. The case of Gabor frames for finite-dimensional
Hilbert spaces was especially treated by
\cite{bava98,bava00,elmawe07,elmasuweXX,feqi95-1,fehakamane07,kust05,li99,pfwa05,qi96,qi98,st97-3}, e.g. the Zibulski-Zeevi
representation of a Gabor frame operator in the finite-dimensional setting was discussed in \cite{elmawe07} and in her Ph.D. thesis \cite{ma07-3} Ewa Matusiak has presented an approach based on representation theory of the Heisenberg group for elementary locally compact abelian groups. Recently,Kaiblinger  has used finite-dimensional Gabor frames to approximate dual Gabor windows of a continuous signal, see \cite{ka05}.

\par
Let $\cH$ be a finite-dimensional Hilbert space of dimension $\rmN$ with
inner product $\langle.,.\rangle_{\mathcal H}$.  A finite family   ${\mathcal G}=\{{\bf g}_{\rm i}\}_{\rm i\in I}$, in $\cH$, with an index set
${\rm I}$ of cardinality ${\rm M}$ is called a {\it frame} if it  is a {\it spanning family}  for $\cH$. Obviously  $\rm M \geq \rm N$.
%
There are three operators naturally associated with  ${\mathcal G}$: the {\it analysis operator} ${\rm C}_{\mathcal G}$, given by
\begin{equation*}
  {\rm C}_{\mathcal G}:  {\bf f}\mapsto (\langle {\bf f},{\bf g}_{\rm i}\rangle_{\mathcal H})_{\rm i\in I}\in\CC^{\rm M},
\end{equation*}
 mapping  $\cH$ into $\rmCCM$ and the {\it synthesis operator}
${\rm D}_{\mathcal G}$ from $\rmCCM$ back  $\cH$, given by
\begin{equation*}
{\rm D}_{\mathcal G: }  {\bf c} =({\rm c_{\rm i}})_{{\rm i\in
I}}\in\rmCCM \mapsto \sum_{\rm i\in I}c_{\rm i}{\bf g}_{\rm i}.
\end{equation*}
An elementary computation shows that
\begin{equation*}
     \langle {\rm C}_{\mathcal G} {\bf f}, {\bf c} \rangle_{\mathcal H} = \langle  {\bf f},  {\rm D}_{\mathcal G} {\bf c} \rangle_{\rmCCM} \,  ,
\end{equation*}
hence these two operators are adjoint to each other,
 ${\rm C}_{\mathcal G}^\star={\rm D}_{\mathcal G}$ and ${\rm D}_{\mathcal G}^\star={\rm C}_{\mathcal G}$.
The most important operator associated to $\{{\bf g}_{\rm i}\}_{\rm i\in I}$
is the {\it frame operator} ${\rm S}_{\mathcal G}$ given  by
\begin{equation*}
  \rmS_{\mathcal G}{\bf f}=\sum_{\rm i\in I}\langle {\bf f},{\bf g}_{\rm i}\rangle_{\mathcal H}\, {\bf g}_{\rm i}.
\end{equation*}
It is easy to show that ${\mathcal G}$ us a frame if and only if $ \rmS_{\mathcal G}$ is invertible. In fact,   the
null-space of $ \rmS_{\mathcal G}= {\rm C}_{\mathcal G}^\star  \circ {\rm C}_{\mathcal G}$ coincides with  the kernel of ${\rm C}_{\mathcal G}$, but  $\operatorname{ker}({\rm C}_{\mathcal G})=\{0\}$ if and only if ${\mathcal G}$ is a frame.
Using next the fact that $\cH$ is finite-dimensional we can argue that the continuous function $ {\bf f} \mapsto \sum_{i\in I}|\langle {\bf f},{\bf
g}_i\rangle_{\mathcal H}|^2  $ is non-zero on the compact
unit-sphere $\{ {\bf f}, \|{\bf f}\| = 1 \}$ of $\cH$, hence
$\{{\bf g}_{\rm i}\}_{\rm i\in I}$ is a frame for $\cH$ if (and only if)
there exist some constants ${\rm A,B}>0$ such that for all non-zero ${\bf
f}\in\cH$ \footnote{In the infinite dimensional setting this is the correct definition.}
\begin{equation*}
 {\rm A\|{\bf f}\|^2_{\mathcal H}\le\|C_{\mathcal G}{\bf f}\|^2_{\rmCCM}=\sum_{i\in I}|\langle {\bf f},{\bf g}_i\rangle_{\mathcal H}|^2 \le B\|{\bf f}\|^2_{\mathcal H}}.
\end{equation*}
If ${\rm A=B}$ in the preceding inequalities, then $\{{\bf g}_{\rm i}\}_{\rm
i\in I}$ is called a {\it tight frame} for $\cH$. Some authors call tight
frames with ${\rm A=B}=1$ {\it normalised tight frames} or {\it Parseval
frames}. Any frame  ${\mathcal G}$  provides natural expansions  of arbitrary
elements ${\bf f}\in\cH$:
\begin{eqnarray*}
  {\bf f}&=&\rmS_{\mathcal G}\rmS_{\mathcal G}^{-1}\,{\bf f}=\sum_{\rm i\in I}\langle {\bf f},\rmS_{\mathcal G}^{-1}{\bf g}_{\rm i}\rangle_{\mathcal H} {\bf g}_{\rm i}\\
  &=&\rmS_{\mathcal G}^{-1}\rmS_{\mathcal G}\,{\bf f}=\sum_{\rm i\in I}\langle {\bf f},{\bf g}_{\rm i}\rangle_{\mathcal H}\, \rmS_{\mathcal G}^{-1}{\bf g}_{\rm i}\\
  &=&\rmS_{\mathcal G}^{-1/2}\rmS_{\mathcal G}\rmS_{\mathcal G}^{-1/2}\,{\bf f}=\sum_{\rm i\in I}\langle {\bf f},\rmS_{\mathcal G}^{-1/2}{\bf g}_{\rm i}\rangle_{\mathcal H}\, \rmS_{\mathcal G}^{-1/2}{\bf g}_{\rm i}.\\
  \end{eqnarray*}
  These formulas also show how the reconstruction of $f$ from the frame coefficients
  $\big(\langle {\bf f},{\bf g}_{\rm i}\rangle_{\mathcal H}\big)_{\rm i\in I} $
   is
  possible.
In general the family $\{{\bf g}_{\rm i}\}_{\rm i\in I}$ is not linear independent which implies the non-uniqueness of the frame decompositions of $\Bff$.
\par
The paper is organized as follows: In Section 2 we recall some basic
facts about finite-dimensional matrix algebras. We emphasize the relevance of the existence of a trace and of the $C^*$-algebra structure of
a finite-dimensional matrix algebra. The material in this  section makes our exposition self-contained and hopefully accessible
to engineers and graduate students, even if they are not yet familiar with classical Gabor analysis. In Section 3  the finite Heisenberg group is introduced, and we present  time-frequency shifts as a projective representation of the time-frequency plane $\rmG\times\widehat\rmG$ (analogue to the Schr\"odinger representation), and its relation  to the twisted group algebra for $\rmG\times\widehat\rmG$. The main result  is Theorem 3.8, describing the  spreading representation for linear operators on $\rmCCN$. In Section 4  we treat (multi-window) Gabor
systems generated with the help of a subgroup $\Lambda\vartriangleleft \rmG\times\widehat \rmG$
and characterize the $\Lambda$-invariant operators by their  spreading support, respectively by their Janssen representations.  We also
discuss the Wexler-Raz biorthogonality relations and the duality principle of
Ron-Shen for finite-dimensional Gabor frames.
Finally, Section 5  contains  results about the optimality of the canonical dual and
tight Gabor window with respect to different measures. In particular we point to the connection
between the tight canonical Gabor window and the L\"owdin orthogonalization, \cite{lo50}.

\section{Basics on matrix algebras}  

Let  $\rmCCN$ and  $\rmCCM$ be endowed with orthonormal
bases $\{\rm e_j:j=1,...,N\}$ and $\{\rm f_j:j=1,...,M\}$. Then the linear mappings $A$
from $\rmCCN$ to  $\rmCCM$ can be
identified with the associated $\rmM\times\rmN$ matrix $\rmA=(\rm a_{\rm i,j})$ with
entries ${\rm a_{\rm i,j}}=\langle\rmA {\rm e_j,f_i}\rangle$ for $\rm
i=1,...,M$ and $j=1,...,N$. We denote with $\rmM\times\rmN$-matrices $\MMNC$ the set of
all those matrices.
The set $\MMC$ is an algebra with respect to matrix multiplication.
with  an {\it involution} $\star$ which for $\rmA=(\rm
a_{\rm i,j})$ is defined by $\rmA^\star=(\ol{\rm a_{\rm j,i}})$. The map
$\rmA\to\rmA^\star$ is an anti-isomorphism, which means that
$(\lambda\rmA)^\star=\ol{\lambda}\rmA^\star$ for $\lambda\in\CC$ and
$(\rmA\rmB)^\star=\rmB^\star\rmA^\star$. Recall that for general $\rmA\in\MMNC$ the
matrices $\rmA\rmA^\star$ and $\rmA^\star\rmA$ are hermitian, and have the
{\it same non-zero eigenvalues} $\lambda_{\rm 1},...,\lambda_{\rm r}$, counting
multiplicities. The square roots of the  eigenvalues $\rmA^\star\rmA$ are
called the {\it singular values} ${\rm s_1\ge s_2\ge\cdots\ge s_r}>0$ where
${\rm r}$ is the rank of $\rmA$. A main result on rectangular matrices $\rmA$
in $\MMNC$ is the {\it singular value decomposition}. It tells us that for
every $\rmA\in\MMNC$ there exist unitary matrices ${\rm U}\in\MMC$ and ${\rm
V}\in\MNC$ and a diagonal matrix ${\rm D}$ with non-negative diagonal entries
$\operatorname{diag}(\rm s_1,...,s_r,0,...,0)$ such that $\rm A={\rm
UDV^\star}$. Despite the possible lack of uniqueness of the factorization,
the sequence of singular values is uniquely determined, because they are
just the eigenvalues of $\rmA^\star\rmA$. We refer the interested reader to \cite{de97} for 
more information about the SVD.
Hence they can be used to define a family of norms $\{\|.\|_{\rm S^p}\}$ on
$\MMNC$,  for $\rm p\in[1,\infty]$, which are indeed the finite-dimensional
versions of the so-called Schatten-von Neumann classes of compact operators
on a Hilbert space. They are given by
 $$ \|\rmA\|_{\rm S^p}=\Big(\sum_{\rm i=1}^{\rm r}{\rm s_i}^{\rm p}\Big)^{\rm 1/p}.$$
The case $\|.\|_{\rm S^2}$ arises in a variety of applications and is
called {\it Frobenius, Hilbert-Schmidt} or {\it Schur norm}. In the following
we will denote $\|\rmA\|_{\rm S^2}$ by $\|\rmA\|_{\operatorname{Fro}}$ and
refer to it as the Frobenius norm of $\rmA\in\MMNC$. The other important
Schatten-von Neumann norm $\|.\|_{\rm S^p}$ arises for $\rm p=\infty$. It turns out to be
just the  {\t operator norm}: 
$$\|\rmA\|_{\operatorname{op}} = \|\rmA\|_{\rm S^\infty}=\lim_{\rm
p\to\infty}\|\rmA\|_{\rm S^p}={\rm s_1} = max \{ s_k: 1 \leq k \leq r \}.$$ 
%
The norms $\|\rmA\|_{\operatorname{Fro}}$ and $\|\rmA\|_{\operatorname{op}}$
  are intimately related to its Euclidean structure and the  $\CS$-algebraic structure of  on  $\MNC$ 
respectively. These two structures will allow us to establish basic facts in a very direct way
(also allowing natural counterparts in the more elaborate case of continuous variables).  
\par
Recall that an {\it algebra} $\A$ is a vector space over $\CC$ with a
multiplication  compatible with the linear structure. If $\A$ has a unit
element, then $\A$ is called {\it unital}. If it is equipped with an {\it involution} ${\rm A}\mapsto{\rm A}^*$ satisfying $({\rm AB})^*={\rm B}^*{\rm A}^*$ and ${\rm A}^{**}={\rm A}$, then $\mathcal{A}$ is called an {\it involutive algebra}. An algebra $\mathcal{A}$ is called a {\it normed algebra}, if there exists a norm $\|.\|_\A$ on $\A$ such that $\|{\rm AB}\|_\A\le \|{\rm A}\|_\A\|{\rm B}\|_\A$. We use  ${\CS}$-algebras, which are involutive normed
algebras $(\A,\|.\|_\A)$ with
$  \|\rmA^\star\rmA\|_\A=\|\rmA\|_\A^2~~\text{for all}~~\rmA\in\A.$
\par
Full matrix algebras $\MNC$ with the operator norm are  $\CS$-algebras, with involution given the transition to transposed and conjugate
 matrices (adjoint operators). Conversely, any
 {\it finite-dimensional $\CS$-algebra} $\A$ is (isometrically) isomorphic to a full matrix algebra
$\MNC$, hence finite-dimensionality forces $\A$ to be unital. We refer the reader to \cite{ta02} for an extensive discussion of finite-dimensional $\CS$-algebras.  
Next we claim  that this are all finite-dimensional $\CS$-algebras. 
This elementary fact will be used  in our approach to Ron-Shen duality
for Gabor frames in $\rmCCN$ below.
\begin{lemma}  The operator norm is the only operator norm on a full matrix algebra turning $\MNC$ into a $\CS$-algebra.
\end{lemma}
\begin{proof}
Since  $\A$ is isomorphic to $\MNC$ for some $\rmN$ our argument relies on the relation between the Schatten-von Neumann classes and the operator norm $\|.\|_{\operatorname{op}}$. Let $\rmA$ be in $\MNC$ and ${\rm D}=\rmA^\star\rmA$. We have to show that $\|{\rm D}\|_\A=\|{\rm D}\|_{\operatorname{op}}$. Let ${\rm d_1\ge d_2\ge\cdots d_r}>0$ be the singular values of $\rmA$, i.e.\  the eigenvalues of $\rmA^\star\rmA$.  The equivalence of any two norms on a finite-dimensional vector space yields:
\begin{equation*}
 {\rm d_1}=\inf_{{\rm R\ge 0}}\Big(\lim_{\rm n\to\infty}\frac{\rm \|D^{2^m}\|_\A}{\rm R^{2^m}}\Big)=\inf_{{\rm R\ge 0}}\Big(\lim_{\rm n\to\infty}\frac{\rm \|D\|_\A^{2^m}}{\rm R^{2^m}}\Big).
\end{equation*}
In the last equality we have used the $\CS$-algebra property of $\|.\|_\A$ and the fact that ${\rm D}={\rm D^2}$ implies $\|{\rm D}^2\|_{\A}=\|{\rm D}\|_\A^2$, hence  ${\rm d_1}=\|{\rm D}\|_\A=\|\rmA^\star\rmA\|_{\operatorname{op}}=\|{\rm
D}\|_{\operatorname{op}}.$
\end{proof}
Next, we recall some well-known facts about traces on matrix algebras. For $\rmA \in\MNC$ the {\it trace} of $\rmA=$ is defined as the sum of its diagonal elements:
\begin{equation*}
  \tr \rmA ={\rm a_{1,1}}+\cdots+a_{\rm N,N}.
\end{equation*}
The trace $\operatorname{tr}$  is a linear functional on $\MNC$ has the following properties:
\begin{enumerate}
  \item $\tr {\rmA^\star}=\ol{\tr \rmA}.$
  \item $\tr{\rmA\rmB}=\tr{\rmB\rmA}$ for all $\rmA,\rmB\in\MNC$ {\it (tracial property)}.
\end{enumerate}
(1) is obvious, and (2) follows from a simple calculation:
\begin{equation*}
  \tr{\rmA^\star\rmB}=\sum_{\rm i,j=1}^{\rm N}\ol{\rm a_{i,j}}{\rm b_{i,j}}=\ol{\tr{\rmB^\star\rmA}}.
\end{equation*}
The identification of $\MNC$ with $\CC^{\rm N^2}$ gives that
$(\rmA,\rmB)\mapsto\tr{\rmA\rmB^\star}$ is an inner product on $\MNC$. In
particular $\operatorname{tr}$ is {\it non-degenerate}, i.e.\
\begin{equation*}
  \tr{\rmA\rmB}=0~~~\text{for all}~~\rmB\in\MNC\Rightarrow\rmA=0.
\end{equation*}
Therefore $(\MNC,\operatorname{tr})$ is a $\rm N^2$-dimensional Hilbert
space. Consequently,  every linear functional
$\phi:\MNC\to\CC$ is of the form $\phi(\rmA)=\tr{\rmA\rmB}$ for some
$\rmB\in\MNC$. The tracial property and Schur's lemma, see Proposition $3.6$  imply
uniqueness of the trace.
\begin{lemma}
Up to multiplicative constants $\operatorname{tr}$ is the only linear functional on $\MNC$ with the tracial property.
\end{lemma}
\begin{proof}
We have to show that for any  linear functional $\phi$ on $\MNC$ with the tracial property $\phi(\rmA\rmB)=\phi(\rmB\rmA)$,  there exists a constant $\rm c\in\CC$ such that $\phi(\rmA)=\rm c\tr\rmA$ for all $\rmA\in\MNC$.
  Let $\rmB$ in $\MNC$ such that $\phi(\rmA)=\tr{\rmA\rmB^\star}$. Then $\tr{(\rm BC-CB)A}=\phi(\rm A C^\star)-\phi(\rm CA^\star)=0$ for all $\rm A,C\in\MNC$. Since trace is non-degenerate $\rm BC=CB$ for all $\rm C\in\MNC$, which means that $\rmB$ is in the center of $\MNC$. Consequently $\rmB=\rm c{\mathbb I}_N$ for some $\rm c\in\CC$.
\end{proof}
Therefore $\langle\rmA,\rmB\rangle_{\Fro}:=\tr{\rmA\rmB^\star}$ is the unique
inner product on $\MNC$ and
\begin{equation*}
  \|\rmA\|_{\Fro}^2 :=\tr{\rmA\rmA^\star}=\tr{\rmA\rmA^\star}=\sum_{\rm i,j=1}^{\rm N}|\rm a_{\rm i,j}|^2
\end{equation*}
defines a norm which  endows $\MNC$ with a {\it unique} Euclidean structure.
\par
For ${\mathcal S}  \subseteq \A$  the {\it subalgebra generated} by the set ${\mathcal S}$ is
denoted by $\operatorname{alg}(\mathcal S)$ and is by definition the smallest
subalgebra of $\A$ that contains ${\mathcal S}$, i.e.\
$\operatorname{alg}(\mathcal S)=\operatorname{span}\{\rm s_1\cdots
s_n:s_1,...,s_n\in{\mathcal S}\}$. Recall also that the {\it commutant}
${\mathcal S}'$ of a set $\mathcal S \subseteq \MNC$ is
\begin{equation*}
  {\mathcal S}' := \{\rmA\in{\mathcal S}:\rm AB=BA~~\text{for all}~~B\in\MNC\}.
\end{equation*}
Observe that ${\mathcal S}'$ is always a subalgebra of $\MNC$ containing all
scalar multiples of the unit in $\MNC$. Moreover the commutant of $\MNC$ is
called the {\it center} of $\MNC$, see \cite{fa01} for more information on
finite-dimensional algebras.
\par
Finally, we want to mention some facts about projective group representations.
Let $(\rmG,\cdot)$ be a finite group and let $\rm V$ be a finite-dimensional Hilbert space.  
Then a {\it projective group representation} $\rho:\rm G\to V$ of $\rmG$ is a family of unitary mappings
$\{\rho(\rm g):\rm g\in\rmG\}$ such that $\rho(\rm g_1\cdot\rm g_2)={\rm
c_G}({\rm g_1,\rm g_2})\rho(\rm g_1)\rho(\rm g_2)$ for unimodular numbers
${\rm c_G}({\rm g_1,\rm g_2})$. The projective group representation $\rho:\rm
G\to V$ of $\rmG$ is called {\it irreducible} if $\{0\}$ and $\rm V$ are the
only $\rho$-invariant subspaces of $\rm V$.
\section{Spreading representation}
In recent years various frame constructions for finite-dimensional Hilbert spaces have been
proposed, using representations of finite Abelian groups (see e.g.\ \cite{boel03}). Because these frames inherit
symmetries from the groups they tend to have interesting extra properties.
We will study frames associated with representations
of the Heisenberg group, which is a non-commutative group, a so-called
two-step nilpotent group, see \cite{grho01,sc04-4} for a detailed treatment
of Heisenberg groups.
\par
Writing  $\TT=\{\tau\in\CC:|\tau|=1\}$ for the circle group,  the {\it
(Weyl)-Heisenberg group} $${\mathbb H}(\rmZZN)=\{(\tau,{\rm k},{\rm
s}):{\rm k},{\rm s}\in\rmZZN,\tau\in\TT\}$$  is endowed with the following multiplication:
\begin{equation}
  (\tau_1,{\rm k}_1,{\rm r}_1)(\tau_2,{\rm k}_2,{\rm r}_2)=(\tau_1\tau_2\cdot e^{2\pi i{\rm k}_2{\rm r}_1},{\rm k}_1+{\rm k}_2,{\rm r}_1+{\rm r}_2).
\end{equation}

Next we fix some notations concerning Gabor frames for
finite-dimensional Hilbert space $\rmCCN$. Its elements are considered as
 discrete time signals of period $\rmN$, i.e.\ ${\rm
f(k)=f(k+mN)}$ for all $\rm m\in{\mathbb Z}$ and if $\rmk$ exceeds the period
$\rmN$ then $\rmk$ is taken modulo $\rmN$. \footnote{So in fact $\rmCCN$ is identified with $\ell^2(\rmZZN)$.}
In the sequel we represent a signal ${\bf f}={\rm
(f(k))}=(f(0),...,f(\rmN-1))^T$ as a {\it column vector} of $\rmCCN$. The
Euclidian structure of $\rmCCN$ induces an inner product on discrete
$\rmN$-periodic signals by $\langle{\bf  f},{\bf g}\rangle_{\rmCCN}=\sum_{\rm
k=0}^{\rm N-1}\rm f(k)\overline{\rm g(k)}$ for ${\bf f},{\bf g}\in\rmCCN$.
\par
The key players of our investigation are {\it time-frequency shifts} of
discrete $\rmN$-periodic signals.
For an integer ${\rm k}$ the {\it translation operator} ${\rm T}_{\rm k}$ is
defined by
\begin{equation*}
  {\rm T}_{\rm k}{\bf  f}=({\rm f({\rm k}),{\rm f}({\rm k}+1),...,{\rm f}({\rm k}-1)}),~~{\bf f}=({\rm f(j)})\in\rmCCN,
\end{equation*}
the {\it modulation operator} $M_{\rm r}$ is given by
\begin{equation*}
  {\rm M}_{\rm r}{\bf  f}=({\rm f(0),e^{2\pi i{\rm r}/\rmN}f(1),e^{2\pi i2{\rm r}/\rmN}f(2),...,e^{2\pi i{\rm r}(\rmN-1)/N}f(N-1)}),~~{\bf f}=({\rm f(j)})\in\rmCCN,
\end{equation*}
and the {\it time-frequency shift} $\pi({\rm k},{\rm r})$ of ${\bf f}$ by
\begin{equation*}
  \pi({\rm k},{\rm r}){\bf  f}={\rm M}_{\rm r}{\rm T}_{\rm k}{\bf  f},~~\text{for}~~{\bf  f}=({\rm f(j)})\in\rmCCN.
\end{equation*}
These operators $\pi({\rm k},{\rm r})$ are called time-frequency shift operators (because time shifts are applied first) and will be the
key-players in our investigation.
\par
We identify $\MNC$ with $\CC^{\rm N^2}$ and denote the standard basis of
$\CC^{\rm N^2}$ by $\{\delta_{\rm i,j}:\rm i,j=1,...,N\}$. Then the matrix
representations of translation $\rm T_1$ and of the modulation ${\rm M_1}$
are given by
$${\rm T_1}=\left(
\begin{array}{ccccc}
0&1&0\cdots&0 \\

\vdots&\ddots&1&\vdots \\
0&0&\cdots& 1 \\
1&0&0\cdots&0
\end{array}\right)~~\text{and}~~{\rm M_1}=\left(
\begin{array}{cccc}
1&0&\cdots&0 \\
0&e^{2\pi i/\rm N}&\ddots&0 \\
\vdots&\ddots&\ddots&\vdots\\
0&0&\cdots&e^{2\pi i{\rm(N-1)/N}}\end{array}\right).$$ An elementary
computation yields ${\rm T_1}^{\rm N}={\rm M_1}^{\rm N}=\mathbb{I}_{\rm N}$
and
\begin{equation*}
  {\rm M_1}{\rm T_1}=e^{2\pi i/\rm N}{\rm T_1}{\rm M_1}.
\end{equation*}
The preceding equation is a finite-dimensional analogue of the
non-commutativity of translation and modulation. In his work on the
foundations of quantum mechanics Weyl treated finite-dimensional analogues of
the commutation relations which appeared in the works of Born, Dirac,
Heisenberg, Jordan and Schr\"odinger on quantum mechanics, \cite{we28}. In
these investigations Weyl was heading towards the fundamental theorem that
all solutions of Born's commutation relation are unitarian equivalent to each
other. He was only able to prove this result for the finite-dimensional case.
Later Stone and von Neumann independently turned Weyl's formal arguments into
a rigorous proof which is the famous Stone-von Neumann theorem. This theorem is one of the most important facts in non-commutative
harmonic analysis and lies at the heart of time-frequency analysis.
\begin{theorem}[Weyl]
Let ${\rm U}, \rm V$ be unitary operators on $\rmCCN$ such that $\text{alg}(\rm U,V)$, the
algebra  generated by $\rm U$ and $ \rm V$ is all of $\MNC$,  and assume they  satisfy the commutation relation
\begin{equation}\label{e:CommutationRelation}
 \rm VU=e^{2\pi i\rm k/N} UV~~\text{for}~~\operatorname{gcd}(k,N)=1.
\end{equation}
Then $\rm U$ and $\rm V$ are unitarily equivalent to $\rm T_1$ and $\rm M_1$,
i.e.\  there exists a unitary operator $Z$ such that ${\rm Z^\star U Z=T_1}$
and ${\rm Z^\star V Z=M_1}$.
\end{theorem}
 Weyl had described the result in terms of projective representations
of $\rmZZN\times\rmZZN$.
\begin{proposition}[Weyl]
Let $ \{\rho({\rm k},{\rm s}):{\rm k},{\rm s}\in\ZZ_{\rm N}\}$ be an
irreducible projective representation of $\rmZZN\times\rmZZN$. Then $\rho$ is
unitarily equivalent to the projective representation of $\rmZZN\times\rmZZN$
by time-frequency shifts $\{\rmpikr={\rm M_rT_k:\rmk,\rmr\in\rmZZN}\}$.
\end{proposition}
It is a fundamental fact of great importance for us that the mapping from $
({\rm k},{\rm r})$ to the time-frequency shift operators $\{\pi({\rm k},{\rm
r}):\rmk,\rmr\in\rmZZN\}$ defines an {\it irreducible projective
representation } of $\rmZZN\times\rmZZN$, i.e.
\begin{equation}
  \pi({\rm k},{\rm r})\pi({\rm l},{\rm s})  \, = \, e^{2\pi i({\rm l\cdot r -k\cdot s})}\pi({\rm l},{\rm s})\pi({\rm k},{\rm r}).
\end{equation}
One advantage of the projective representation instead of unitary
representation for the Heisenberg group of $\rmZZN$ is that the
matrix-coefficients $\langle\Bff,\pi({\rmk,\rmr})\Bfg\rangle$ for
$\Bff,\Bfg\in\rmCCN$ of $\{\pi(\rmk,\rmr):\rml,\rmr\in\rmZZN\}$ have a
concrete practical meaning and we do not have to take care of additional
phase factors. Namely as the {\it Short-time Fourier transform} of
$\Bff\in\rmCCN$ with respect to a fixed {\it window} $\Bfg\in\rmCCN$
\begin{equation*}
  {\rm V}_{\Bfg}\Bff(\rmk,\rmr)=\sum_{\rml\in\rmZZN}{\rm f}(\rm l)\overline{g}(\rm k-l)e^{-2\pi i\rmr\rml/N}=\langle\Bff,\pi(\rmk,\rmr)\Bfg\rangle_{\rmCCN}.
\end{equation*}
An important fact about STFT's is the so-called {Moyal Identity}.
\begin{proposition}[Moyal's Formula]\label{moyal}
Let $\Bff_1,\Bfg_1$ and $\Bff_2,\Bfg_2$ in $\rmCCG$. Then
\begin{equation*}
  \langle{\rm V}_{\Bfg_1}\Bff_1,{\rm V}_{\Bfg_2}\Bff_2\rangle_{\rmCCGG}=\langle\Bff_1,\Bff_2\rangle_{\rmCCG}\overline{\langle\Bfg_1,\Bfg_2\rangle_{\rmCCG}}.
\end{equation*}
\end{proposition}
The proof is straightforward and relies on the definition of the STFT and a
change of summation. Technically speaking Moyal's identity expresses the
orthogonality of two matrix coefficients of the irreducible projective
representation $\{\ulrmpikr:(\ulrmk,\ulrmr)\in\rmGG\}$ of $\rmGG$. An
application of Moyal's Identity with $\Bff=\Bff_1=\Bff_1$ and
$\Bfg=\Bfg_1=\Bfg_2$ yields that ${\rm C}=\|\Bfg\|_{\rmCCG}^2$. 
\par
If we regard the elements of $\rmZZN\times\rmZZN$ as linearly independent
vectors, most conveniently as $\delta_{\rm k,0}\cdot\delta_{\rm r,0}$ for
${\rm k,r\in\rmZZN}$, then their span
is the vector space
\begin{equation*}
  \A(\rmZZN\times\rmZZN)=\Big\{\sum_{\rmk,\rmr\in\rmZZN}\rmakr\delta_{\rm k,0}\cdot\delta_{\rm r,0}\Big\}.
\end{equation*}
The vector space $\A(\rmZZN\times\rmZZN)$ has a product which arises from the
group product by a twist of a uni-modular number. Therefore the arising group
algebra is referred to as the {\it twisted group algebra} of
$\rmZZN\times\rmZZN$. More precisely, the multiplication of $\delta_{\rm
k,0}\delta_{\rm r,0}$ and $\delta_{\rm l,0}\delta_{\rm s,0}$ is the {\it
(left) twisted translation} $\delta_{\rm k,0}\delta_{\rm r,0}\delta_{\rm
l,0}\delta_{\rm s,0}=e^{2\pi i\rm (l-k)\cdot r/N}\delta_{\rm
k+l,0}\delta_{\rm r+s,0}$. For the elements of  $\A(\rmZZN\times\rmZZN)$ this
yields
\begin{eqnarray*}
  &&\Big(\sum_{\rm k,r}\rmakr\delta_{\rm k,0}\cdot\delta_{\rm r,0}\Big)\Big(\sum_{\rm l,s}\rmbls\delta_{\rm l,0}\cdot\delta_{\rm s,0}\Big)\\
  &=&\sum_{\rm k,r}\sum_{\rm l,s}\rmakr\rmbls e^{2\pi i\rm r\cdot l/rmN}\delta_{\rm k+l,0}\cdot\delta_{\rm r+s,0}\\
  &=&\sum_{\rm k,r}\Big(\sum_{\rm l,s}{\rm a}(\rm l,\rm s){\rm b(l-k,s-r)}e^{2\pi i\rm (l-k)\cdot r/N}\delta_{\rm k,0}\cdot\delta_{\rm r,0}\Big).
\end{eqnarray*}
This  calculation motivates the definition of a "twisted" product for
$\A(\rmZZN\times\rmZZN)$:
\begin{definition}
The complex vector space $\A(\rmZZN\times\rmZZN)$ of functions on the group 
$\rmZZN\times\rmZZN$ is called the twisted group algebra when given the {\it
twisted convolution} $\natural$ as a product. For ${\bf a}=(\rmakr)$ and ${\bf
b}=(\rmbls)$ in $\A(\rmZZN\times\rmZZN)$ we define  the {\bf twisted
convolution} of ${\bf a}$ and ${\bf b}$ by
\begin{equation}
  ({\bf a}\natural{\bf b})(\rm k,r)=\sum_{\rm l,s}\rmals{\rm b(l-k,s-r)}e^{2\pi i\rm (l-k)\cdot s}.
\end{equation}
\end{definition}
A projective representation of $\rmZZN\times\rmZZN$ induces a representation
of the twisted group algebra ${\mathcal A}(\rmZZN\times\rmZZN)$. We define a
representation $\pi_\A$ for ${\mathcal A}(\rmZZN\times\rmZZN)$ by
\begin{equation*}
  \pi_\A(\bf a): = \sum_{{\rm k},{\rm r}\in\rmZZN}{\rm a}({\rm k},{\rm r})\pi({\rm k},{\rm r}) \quad \text{for} \, {\bf a}\in{\mathcal A}(\rmZZN\times\rmZZN).
\end{equation*}
For non-discrete groups the above representation of the twisted group algebra
is the so-called {\it integrated representation}. The  
representation $\pi_\A$ of ${\mathcal A}(\rmZZN\times\rmZZN)$ is intimately
related with the representation of the group $\rmZZN\times\rmZZN$. Therefore
we collect some of the basic properties of time-frequency shifts which are
elementary consequences of the commutation relations.
\begin{lemma}\label{TFprop}
Let $(\rm k,r)$ and $(\rm l,s)$ be in $\rmZZN\times\rmZZN$. Then
\begin{enumerate}
  \item[(i)]  $\pi(\rm k,r)^*={\rm T_{-\rm k}}{\rm M_{-\rm r}}=e^{2\pi i\rm k\cdot r/N}\pi(\rm -k,-r),$
  \item[(ii)] $\pi(\rm k,r)\pi(\rm l,s)=e^{2\pi i\rm l\cdot r/N}\pi(\rm k+l,r+s),$
  \item[(iii)] $\pi(\rm k,r)\pi(\rm l,s)=e^{2\pi i(\rm l\cdot r-\rm k\cdot s)/N}\pi(\rm l,s)\pi(\rm k,r).$
\end{enumerate}
\end{lemma}
Lemma \ref{TFprop} (i) suggests the following natural {\it twisted
involution}: 
on ${\mathcal A}(\rmZZN\times\rmZZN)$.  
\begin{equation}\label{twistinvol}
  {\bf a}^\star(\rm k,r) = e^{2\pi i\rm k\cdot r/N} \, \ol{\rm a(-k,-r)}.
\end{equation}

Altogether we have  the following properties of $\pi_\A$:
\begin{proposition} The mapping $ {\bf a} \mapsto \pi_\A(\bf a)$, defined on
$\A(\rmZZN\times\rmZZN)$, given by
\begin{equation}\label{invrep}
 \pi_\A(\bf a)=\sum_{{\rm k},{\rm r}\in\rmZZN}{\rm a}({\rm k},{\rm r})\pi({\rm k},{\rm r})
\end{equation}
defines an involutive representation of $\A(\rmZZN\times\rmZZN)$, i.e.\ one
has:
\begin{enumerate}
  \item[(a)] $\pi_\A({\bf a}) + \pi_\A({\bf b})=\pi_\A({\bf a}+{\bf b}),$
  \item[(b)] $\pi_\A({\bf a}) \, \pi_\A({\bf b})=\pi_\A({\bf a}\natural{\bf b}),$
  \item[(c)] $\pi_\A({\bf a}^\star)=\pi_\A({\bf a})^\star,$
  \item[(d)] $\pi_\A({\delta_{0,0}})={\rm 1}_{\rm N}.$
\end{enumerate}
Conversely, if an involutive representation $\pi_\A$ of
$\A(\rmZZN\times\rmZZN)$ satisfies $(a)-(d)$, then there exists a 
projective representation $\pi$ of $\rmZZN\times\rmZZN$ satisfying
\eqref{invrep}, i.e., all such algebra representations arise in this way.
\end{proposition}
We still have to  show the converse,   
for which we need the definition of a projective representation $\pi$ for
$\rmZZN\times\rmZZN$ by \eqref{invrep}.
\begin{equation}\label{def-repr2}
  \pi(\rm k,r):=\pi_\A(\delta_{\rm k,\, 0}\cdot\delta_{\rm r,\, 0}).
\end{equation}
%
The situation just described is a special case of a more general result,
according to which there is a one-to-one correspondence between projective
representations of a finite group $\rmG$ and involutive representations of
its twisted group algebra $\A(\rmG\times\rmG)$. Another close relation
between a group and its group algebra is the fact that $\pi_\A$ is an
irreducible involutive representation of $\A(\rmG\times\rmG)$ if and only if
$\pi$ is an irreducible projective representation of $\rmG\times\rmG$.
Schur's lemma formulates this more precisely:
\begin{proposition}[Schur's Lemma]
Let $\pi_\A$ be an irreducible involutive representation of
$\A(\rmZZN\times\rmZZN)$ on $\rmCCN$. If a linear mapping $\rmA$ of $\rmCCN$
satisfies
\begin{equation*}
  \pi_\A({\bf a}) \, \rmA = \rmA  \, \pi_\A({\bf a}) \quad \text{for all}~~{\bf a}\in\A(\rmZZN\times\rmZZN)
\end{equation*}
or equivalently $\pi(\rm k,r)\rmA=\rmA\pi(\rm k,r)$ for all $\rm
k,r\in\rmZZN$. Then $\rm A=\rm c{\mathbb I}_N$ for some $\rm c\in\CC$.
\end{proposition}
For the proof we refer the reader to the excellent book \cite{te99} by
A.~Terras, which provides further information about representations of
finite groups.

By construction the twisted group algebra $\A(\rmZZN\times\rmZZN)$ coincides
with the full matrix algebra $\MNC$.  From this point of view Schur's lemma
states that the center of $\MNC$ consists of all scalar multiples of the
identity $\{\rm c {\mathbb I}_N :c\in\CC\}$.
\par
In our discussion of the spreading function we make use of the existence of a
{\it trace} on $\MNC$. The Hilbert space $\MNC$ with
$\langle\rmA,\rmB\rangle_{\Fro}=\operatorname{tr}(\rmA\rmB^\star)$ has a
natural orthonormal basis which leads to the spreading representation of
linear operators on $\rmCCN$.
\begin{lemma}
Let $\{\pi({\rm k},{\rm r}):{\rm k},{\rm r}\in\rmZZN\}$ be the family of all
time-frequency shift operators in $\MNC$. Then $\{{\rm N}^{\rm -1/2}\pi({\rm
k},{\rm r}):{\rm k},{\rm r}\in\rmZZN\}$ is an orthonormal basis for
$(\MNC,\|.\|_{\rm Fro})$.
\end{lemma}
\begin{proof}
The observation that the cardinality of $\{\pi({\rm k},{\rm r}):{\rm k},{\rm
r}\in\rmZZN\}$ is equal to the dimension of $\MNC$ and the following
calculation gives the desired assertion. Let $\pi({\rm k},{\rm r})$ and
$\pi({\rm l},{\rm s})$ be two time-frequency shifts for ${\rm k},{\rm
r}\in\rmZZN$ and $({\rm l},{\rm s})\in\rmZZN$. Then
\begin{eqnarray*}
 \langle\pi({\rm k},{\rm r}), \pi({\rm l},{\rm s})\rangle&=& \operatorname{tr}(\pi({
 \rm k},{\rm r})\pi({\rm l},{\rm s})^*)\\
 &=&\operatorname{tr}(\pi({ \rm k}-{\rm l},{\rm r}-{\rm s}))={\rm N} \delta_{k-l,s-r}.
\end{eqnarray*}
\vspace{-5mm}
\end{proof}
As a consequence we are able to expand every operator $\rmA\in\MNC$ with
respect to {\it all} time-frequency shifts from the discrete time-frequency
plane $\rmZZN\times\rmZZN$.
\begin{theorem}[Spreading representation]
For $\rmA\in\MNC$ we have
\begin{eqnarray*}
  \rmA&=&\sum_{{\rm k},{\rm r}\in\rmZZN}\langle \rmA,\pi({\rm k},{\rm r})\rangle_{\Fro}\pi({\rm k},{\rm r})\\
   &=&\sum_{{\rm k},{\rm r}\in\rmZZN}\eta_{\rmA}({\rm k},{\rm r})\pi({\rm k},{\rm r}),
\end{eqnarray*}
where $\eta_{\rmA}=(\eta_{\rmA}({\rm k},{\rm r}))_{{\rm k},{\rm r}\in\rmZZN}$
is called the \textbf{spreading function} of $\rmA$. Furthermore
\begin{equation}
\eta_{\rmA} ({\rm k},{\rm r})={\rm N}^{-1}\sum_{{\rm l}\in\rmZZN}{\rm a}({\rm
l},{\rm l}-{\rm k})e^{-2\pi i{\rm l r}/\rm N} \quad \mbox{for} \quad {\rm
k},{\rm r}\in\rmZZN.
\end{equation}

\end{theorem}
The expression for the spreading function is a direct consequence of
the definitions. Traditionally  a linear operator $\rmA$ on $\rmCCN$ is
formulated via a kernel by
\begin{equation*}
  \rmA {\rm f}({\rm j})=\sum_{\rm i\in\rmZZN}{\rm k_A(i,j)}{\rm f(i)},~~\text{for}~~{\bf f}=({\rm f(j)})\in\rmCCN.
\end{equation*}
Then the relation between the spreading representation and the kernel of
$\rmA$ is
\begin{equation*}
  \eta_{\rmA}({\rm k,r})=\sum_{\rm i\in\rmZZN}{\rm k_A}({\rm i,i-k})e^{-2\pi i{\rm r\cdot i/N}}
\end{equation*}
and the inversion formula is
\begin{equation*}
  {\rm k_A}({\rm i,j})=\sum_{{\rm r\in\rmZZN}}\eta_{\rmA}({\rm i-j,r})e^{2\pi i{\rm r\cdot i/N}}.
\end{equation*}
The correspondence between $\rmA\in\MNC$ and its spreading coefficients
$(\eta_{\rmA}(\rmk,\rmr))$ may be considered as non-commutative analogue of
the finite Fourier transform. First we have as a non-commutative analogue of
Parseval's Theorem which follows from the orthogonality of
$\{\rmpikr:\rmk,\rmr\in\rmZZN\}$ and $\tr{A}=\rmN\eta_{\rmA}(0,0)$:
\begin{equation*}
  \langle \rmA,\rmA\rangle_{\Fro}=\tr{\rmA^\star\rmA}=\rmN\sum_{\rm k,r\in\rmZZN}|\eta_{\rmA}(\rmk,\rmr)|^2~~\text{for}~~\rmA\in\MNC.
\end{equation*}
  Recall that the Fourier coefficients provide the best least square approximation under all trigonometric polynomials.
  An analogous argument yields the same for general orthonormal bases in
  abstract Hilbert space. For the spreading function this implies:
\begin{theorem}
Let $\rmA$ be in $\MNC$. Then 
for every  subset ${\rm F} \subseteq \rmZZN \times \rmZZN$ the best
approximation to $\rmA$   among all finite linear combinations of
time-frequency shifts from $F$  - in the Frobenius norm  $\|.\|_{\Fro}$  - is
given by
  $ \sum_{{\rm F}\subseteq\rmZZN}\rmakr\pi({\rm k},{\rm r}).$ 
\end{theorem}
\begin{proof}
The orthogonality of the time-frequency shifts $\{\rmpikr:\rm k,r\in\rmZZN\}$
implies that it is sufficient to prove the following: Let $\rm U$ be a
unitary matrix in $\MNC$. Then $\rmA$ has to be approximated by a scalar
multiple of the unitary matrix $\rm U$.
\begin{equation*}
  \|\rmA-\rm cU\|_{\Fro}^2=\tr{(\rm A-cU)^*(\rm A-cU)}=\|\rmA\|_{\Fro}^2-2\operatorname{Re}{\rm\ol{c}\tr{\rm AU^\star}}+\rm N|\rm c|^2
\end{equation*}
is minimized for $\rm c={\rm N}^{-1}\tr{\rm AU^\star}$. The choice $\rm
U=\rmpikr$ for some $\rm k,r\in\rmZZN$ gives the desired assertion since
$\eta_{\rmA}(\rm k,r)={\rm N}^{-1}\tr{\rm A\pi(\rm k,r)^\star}$.
\end{proof}
The spreading representation  was introduced by Kailath in the electrical
engineering context of time-variant systems \cite{ka62} completely
independent from mathematical physics and representation theory. The
parallelism with representation theory of the Heisenberg group becomes
evident through the work of Feichtinger and Kozek in \cite{feko98}. A real
world realization of time-variant filters requires a finite-dimensional
model which can be implemented on a computer. The spreading representation
allows a symbolic calculus for time-variant filters which in mathematics are
called pseudo-differential operators.
\par
Let $\rmA$ be our finite-dimensional model of a time-variant system. Then the
spreading function $\eta_{\rmA}$ can be considered as a symbol which contains
all information about the time-frequency concentration of $\rmA$.
Given the spreading representations of $A=\sum_{{\rm k},{\rm
r}\in\rmZZN}\eta_{\rmA}({\rm k},{\rm r})\pi({\rm k},{\rm r})$ and
$B=\sum_{{\rm k},{\rm r}\in\rmZZN}\eta_{\rmB}({\rm k},{\rm r})\pi({\rm
k},{\rm r})$, then after an easy direct computation  of the spreading
representation $AB$ reveals that \begin{equation} \label{comptwist1}
  \rmA\rmB=\sum_{{\rm k},{\rm
r}\in\rmZZN}(\eta_{\rmA}\natural\eta_{\rmB})({\rm k},{\rm r})\pi({\rm k},{\rm
r}).
\end{equation}

These observations allow us to complement recent work of Wildberger on a
symbolic calculus for finite abelian groups. In \cite{wi05} the Weyl quantization
was considered on finite abelian groups and one of the main
results states that there is no good symbolic calculus for groups of even
order. Our results do not rely on the order of $\rmZZN$ which is possible
because we use a different kind of quantization, namely the  Kohn-Nirenberg
quantization. In \cite{grstXX} the Kohn-Nirenberg quantization is discussed
in full generality for locally compact abelian groups, and its relevance for
time-frequency analysis is shown in \cite{st06}. The presentation given here
makes use of the finite-dimensionality of the space of time-variant filters,
i.e.\ linear operators on $\rmCCN$.
\par
The singular value decomposition of $\rmA$ of rank $\rmr$ in $\MNC$ may be
considered as the decomposition of $\rmA$ as the sum of $\rmr$ rank-one
operators ${\bf g}_{\rm i}\otimes\ol{{\bf h}_{\rm i}}$ for ${\bf g}_{\rm
i},{\bf h}_{\rm i}\in\rmCCN$ for ${\rm i}=1,...,{\rm r}$ where ${\bf
g}\otimes\ol{{\bf h}}$ denotes the rank-one operator ${\rm P}_{{\bf g},{\bf
h}}$
\begin{equation}\label{rank-onex}
{\rm P}_{{\bf g},{\bf h}}{\bf f}:=\big({\bf g}\otimes\ol{{\bf h}}\big){\bf
f}=\langle{\bf f},{\bf h}\rangle{\bf g}~~\text{for}~~{\bf h}\in\rmCCN.
\end{equation}
More concretely, if $\rmA$ has the singular value decomposition ${\rm UDV^*}$
with singular values ${\rm d_1\ge d_2\ge\cdots\ge d_r}>0$, then ${\bf g}_{\rm
i}$ and $\ol{{\bf h}_{\rm i}}$ are the ${\rm i}$-th column ${\bf u}_{\rm i}$
of ${\rm U}$ and ${\rm i}$-th row ${\bf v}_{\rm i}$ of ${\rm V}$, i.e.\
\begin{equation*}\label{e:FiRaOn}
  \rmA={\rm d_1}{\bf u}_{\rm 1}\otimes\ol{{\bf v}_{\rm 1}}+\cdots+{\rm d_r}{\bf u}_{\rm r}\otimes\ol{{\bf v}_{\rm r}}.
\end{equation*}
It thus turns out that the spreading representation of $\rmA$ is given by
\begin{equation*}
  \rmA=\sum_{\rm i=1}^{\rmr}\sum_{\rml,\rms\in\rmZZN}{\rm d_i}\langle{\rm P}_{{\bf u}_{\rm i},{\bf v}_{\rm i}},\rmpils\rangle_{\Fro}\rmpils,
\end{equation*}
i.e.\  the problem is reduced to the spreading representation of a rank-one
operator ${\rm P}_{{\bf g},{\bf h}}$ for ${\bf g}$ and ${\bf h}$ as mentioned
above. By definition we get
\begin{equation}\label{spread-repr1}
  {\rm P}_{{\bf g},{\bf h}}=\rmN^{-1}\sum_{\rm k,r\in\rmZZN}\langle {\bf h},\rmpikr {\bf g}\rangle_{\Fro} \, \rmpikr.
\end{equation}
The symbolic calculus for linear operators allows us to transfer properties of
the operators with relations for their spreading functions. As an example we
look at the product of two rank-one operators and their spreading
representation. Let ${\bf g}_{\rm 1},{\bf g}_{\rm 2},{\bf h}_{\rm 1},{\bf
h}_{\rm 2}$ be in $\rmCCN$. Then an elementary computation gives that
\begin{equation}\label{spread-repr2}
  {\rm P}_{{\bf g}_{\rm 1},{\bf h}_{\rm 1}}{\rm P}_{{\bf g}_{\rm 2},{\bf h}_{\rm 2}}=
  \langle{\bf g}_{\rm 2},{\bf h}_{\rm 1}\rangle  {\rm P}_{{\bf g}_{\rm 1},{\bf h}_{\rm 2}}.
\end{equation}
The spreading representation of a rank-one operator and the symbolic calculus
for linear operators yields the so-called {\it reproducing property}
\begin{equation*}
  \langle{\bf g}_{\rm 1},\pi(k,r){\bf h}_{\rm 1}\rangle\natural\langle{\bf g}_{\rm 2},\pi(k,r){\bf h}_{\rm 2}\rangle=\langle{\bf g}_{\rm 2},{\bf h}_{\rm 1}\rangle\langle{\bf g}_{\rm 1},\pi(k,r){\bf h}_{\rm 2}\rangle.
\end{equation*}
In our discussion of Gabor frames we will return to the coefficients of the
spreading representation of a rank-one operator and place it into the setting of
representation theory for a pair of twisted group algebras. At the moment we
want to emphasize the role of projection operators and their associated
one-dimensional subspaces. Our discussion indicates the connection between
our presentation of the spreading representation and the work of Calderbank
in \cite{cahomo03} and Grassmannian frames \cite{hest03}.
\subsection{Discrete time-frequency plane $\rmG\times\widehat \rmG$}
So far we have restricted ourselves to cyclic groups $\rmZZN$. This is no
real restriction, since the structure theorem for finite abelian groups
$\rmG$ allows us to move on to the more general setting of a finite abelian
group $\rmG$.
\begin{lemma}[Decomposition of Finite Abelian Groups]
Let $\rmG$ be a finite Abelian group of order $\rmN$. Then $\rmG$ is
isomorphic to a direct product of groups $\ZZ_{\rm  p_i^{\rm N_i}}$ where
${\rm N=p_1^{N_1}\cdots p_k^{N_k}}$ is the prime number decomposition of $N$,
i.e.\  $\rmG\cong \ZZ_{\rm p_1^{N_1}}\times\cdots\times\ZZ_{\rm p_k^{N_k}}$.
\end{lemma}
A character ${\rm c_G}$ of a finite abelian group is a mapping from
$\rmG\to\TT$,  such that ${\rm {\rm c}_G(x+y,\om)={\rm c}_G(x,\om){\rm
c}_G(y,\om)}$ for all ${\rm x,y\in G}$. The set of all characters of $\rmG$
form a group with respect to multiplication which we denote by
$\widehat\rmG$. A basic result about finite abelian groups asserts that
$\widehat\rmG$ is naturally isomorphic to $\rmG$. In addition we have an
explicit knowledge of ${\rm c_G}$ in terms of the building blocks of $\rmG$,
i.e.\
\begin{equation*}
{\rm c_G(x,\om)={\rm c}_{\ZZ_1}({\rm l}_1,{\rm s}_1)\cdots {\rm
c}_{\ZZ_k}({\rm l}_k,{\rm s}_k)~~x=({\rm l}_1,\cdots,{\rm l}_k)\in
G,\om=({\rm s}_1,\cdots,{\rm s}_k)\in\widehat G},
\end{equation*}
where we abbreviated $\ZZ_{\rm p_i^{N_i}}$ by $\ZZ_{\rm i}$.
\par
We refer to $\rmGG$ as the discrete time-frequency plane. We are therefore
faced with the discussion of twisted group algebras $\A(\rmG_1\times\rmG_2)$
for $\rmG_1$ for abelian groups $\rmG_1$ and $\rmG_2$ which turns out to be
the tensor product $A(\rmG_1\times\rmG_1)\otimes\A(\rmG_2\times\rmG_2)$.
Consequently the twisted group algebra $\A(\rmGG)=\A(\ZZ_{\rm
p_1^{N_1}}\times\ZZ_{\rm p_1^{N_1}})\otimes\cdots\otimes\A(\ZZ_{\rm
p_k^{N_k}}\times\ZZ_{\rm p_k^{N_k}})$ is isomorphic to a product of
full matrix algebras ${\mathcal M}_{|\rmG\times \widehat \rmG|}={\mathcal
M}_{\rm p_1^{N_1}\times p_1^{N_1}}\otimes\cdots\otimes{\mathcal M}_{\rm
p_k^{N_k}\times p_k^{N_k}}$. In the following we denote an element of $\rmGG$
by $(\ulrmk,\ulrmr)$ to emphasize that $\ulrmk,\ulrmr$ are elements of
$\rmG\cong \ZZ_{\rm p_1^{N_1}}\times\cdots\times\ZZ_{\rm p_k^{N_k}}$.
\section{Gabor frames}
In this section we discuss Gabor frames from a matrix algebra point of view.
If $\mathcal{G}(\Bfg,\Lambda)$ is a Gabor frame for $\rmCCG$, then the Gabor
system $\mathcal{G}(g,\Lambda^\circ)$ over the adjoint subgroup
$\Lambda^\circ$ is naturally related to the original Gabor system. We will
consider the twisted group algebras of $\Lambda$ and $\Lambda^\circ$ and the
investigation of their structure allows a unified treatment of Gabor frames.
The understanding of the twisted group algebras for the subgroups $\Lambda$
and $\Lambda^\circ$ requires some harmonic analysis over the time-frequency
plane. We therefore develop this prerequisites in the first part of this
section. In the second part we introduce multi-window Gabor frames and as a
motivation we treat the set of all time-frequency shifts as the Gabor system
${\mathcal G}(\Bfg,\rmGG)$ and state the {\it Resolution of Identity} as a
"continuous" analogue for reconstruction formulas for a general Gabor system.


\subsection{Harmonic Analysis over $\rmGG$}
Harmonic analysis over the time-frequency plane $\rmGG$ is the study of
the Fourier transform of $\rmGG$ and its properties which differ from the
Euclidean Fourier transform of $\rmG\times\rmG$. The difference between
$\rmG\times\rmG$ and $\rmGG$ arises from the symplectic structure of the
time-frequency plane. We motivate our investigation of $\rmGG$ with the
commutation relations for time-frequency shifts:
\begin{equation*} \ulrmpils[({\ulrmk},{\ulrmr})]={\rm c}_{\rmG}\big(({\ulrmk},{\ulrmr}),({\ulrml},{\ulrms})\big)\overline{{\rm c}}_{\rmG}\big(({\ulrml},{\ulrms}),({\ulrmk},{\ulrmr})\big)\ulrmpils=e^{2\pi i\Omega(({\ulrmk},{\ulrmr}),({\ulrml},{\ulrms}))/|\rmG|}\ulrmpils,
\end{equation*}
where $\Omega:\rmGG\to\rmGG$ denotes the {\it standard symplectic form}
\begin{equation*}
  \Omega(({\ulrmk},{\ulrmr}),({\ulrml},{\ulrms})):=\ulrml\cdot\ulrmr-\ulrmk\cdot\ulrms.
\end{equation*}
Therefore the commutation relations for time-frequency shifts endow the
discrete time-frequency plane $\rmGG$ with an intrinsic {\it symplectic
structure}.
\par
From the algebraic point of view $\rmGG$ is isomorphic to $\rmG\times\rmG$
but from a harmonic analysis point of view these two groups are different
objects. Namely, the characters of $\rmGG$ are not of the form ${\rm
c}_{\rmG\times\rmG}={\rm c}_{\rmG}\cdot{\rm c}_{\rmG}$, i.e. the character
group of $\rmGG$ is not $\widehat\rmG\times\widehat\rmG$. Observe that in the
time-frequency plane the character group $\widehat \rmG$ and $\rmG$ are
orthogonal to each other, i.e. a rotation by $\pi/2$ moves $\rmG$ onto
$\widehat\rmG$. These facts should convince you that the "correct" characters
of $\rmGG$ are Euclidean characters rotated by $\pi/2$. More precisely, we
define a {\it symplectic character} ${\rm c}_{|\rmGG|}^{\bf
s}\big(({\ulrmk},{\ulrmr}),({\ulrml},{\ulrms})\big)$ for a fixed
$({\ulrml},{\ulrms})$ by
\begin{equation*}
  {\rm c}_{|\rmGG|}^{\bf s}\big(({\ulrmk},{\ulrmr}),({\ulrml},{\ulrms})\big)={\rm c}_{\rmG}\big(({\ulrmk},{\ulrmr}),({\ulrml},{\ulrms})\big)\overline{{\rm c}}_{\rmG}\big(({\ulrml},{\ulrms}),({\ulrmk},{\ulrmr})\big).
\end{equation*}
\begin{lemma}
The character group of $\rmGG$ is $\{{\rm c}_{|\rmGG|}^{\bf
s}\big(({\ulrmk},{\ulrmr}),({\ulrml},{\ulrms})\big):({\ulrmk},{\ulrmr})\in\rmGG\}$,
i.e. it is isomorphic to $\widehat{\rmG}\times\rmG$.
\end{lemma}

This observation motivates the following "symplectic" analogues of
translation and modulation operators which are actually the translation and
modulation operators for $\rmGG$. Let ${\bf F}$ be in $\rmCCGG$. Then we
define the {\it symplectic translation} operator by
\begin{equation*}
   {\rm T}^{\bf s}_{({\ulrmk'},{\ulrmr'})}{\bf F}({\ulrmk},{\ulrmr})={\bf F}(({\ulrmk}-{\ulrmk'},{\ulrmr}-{\ulrmr'}))~~\text{for}~~({\ulrmk'},{\ulrmr'})
\end{equation*}
and the {\it symplectic modulation} operator by
\begin{equation*}
  {\rm M}^{\bf s}_{({\ulrml},{\ulrms})}{\bf F}(({\ulrmk},{\ulrmr}))=e^{2\pi i\Omega(({\ulrmk},{\ulrmr}),({\ulrml},{\ulrms}))/|\rmG|}{\bf F}\big(({\ulrmk},{\ulrmr})\big)
\end{equation*}
for $({\ulrmk},{\ulrmr}),({\ulrml},{\ulrms})\in\rmGG$. As for translation and
modulation operators we have the non-commutativity of symplectic translations
and modulations:
\begin{equation*}
    {\rm T}^{\bf s}_{({\ulrmk},{\ulrmr})}{\rm M}^{\bf s}_{({\ulrml},{\ulrms})}={\rm c}_{|\rmGG|}^{\bf s}\big(({\ulrmk},{\ulrmr}),({\ulrml},{\ulrms})\big){\rm M}^{\bf s}_{({\ulrml},{\ulrms})}{\rm T}^{\bf s}_{({\ulrmk},{\ulrmr})}
\end{equation*}
In other words $\{{\rm T}^{\bf s}_{({\ulrmk},{\ulrmr})}{\rm M}^{\bf
s}_{({\ulrml},{\ulrms})}:({\ulrmk},{\ulrmr})({\ulrml},{\ulrms})\in\rmGG\}$ is
a projective representation of $\rmGG\times\rmGG$. Therefore we may define
the symplectic analogue of the short-time Fourier transform and analyze
objects on $\rmGG$ with respect to this symplectic short-time Fourier
transform. Gr\"ochenig uses this object implicitly in his discussion of
time-frequency localization operators and pseudo-differential operators,
\cite{cogr03}.
\par
The structure of the character group of $\rmGG$ indicates that the "correct"
Fourier transform in time-frequency analysis is not the standard Euclidean
Fourier transform but the {\it symplectic Fourier transform} of ${\bf
F}\in\rmCCGG$:
\begin{equation*}
  {\mathcal F}_{\bf s}({\ulrmk},{\ulrmr})=|\rmGG|^{-1/2}\sum_{({\ulrml},{\ulrms})}{\rm c}_{|\rmGG|}^{\bf s}\big(({\ulrmk},{\ulrmr}),({\ulrml},{\ulrms})\big){\bf F}\big(({\ulrml},{\ulrms})\big).
\end{equation*}
An elementary computation establishes the following property of the
symplectic Fourier transform.
\begin{lemma}
The symplectic Fourier transform ${\mathcal F}_{\bf s}$ is a self-inverse
mapping of order two on $\rmCCGG$, i.e.
\begin{equation*}
  {\mathcal F}_{\bf s}^{-1}={\mathcal F}_{\bf s},\qquad{\mathcal F}_{\bf s}^2={\rm I}_{\rmCCGG
}.
\end{equation*}
\end{lemma}
Since the symplectic characters ${\rm c}_{|\rmGG|}^{\bf s}$ are obtained from
the characters ${\rm c}_{|\rmG\times\rmG|}$ by a rotation of $\pi/2$ the
Euclidean Fourier transform ${\mathcal F}$ and the symplectic Fourier
transform ${\mathcal F}_{\bf s}$ are related by
\begin{equation*}
{\mathcal F}_{\bf s}={\mathcal F}\circ{\rm J}~~\text{for}~~  {\rm J}=\left(
\begin{array}{cc}
0&{\rm I}_{|\rmG|}\\
-{\rm I}_{|\rmG|}&0
\end{array}\right),
\end{equation*}
where ${\rm J}$ describes the rotation by $\pi/2$ of the time-frequency plane
$\rmGG$. We want to have a Poisson summation formula for the symplectic
Fourier transform but what is the symplectic analogue of the dual subgroup
$\Lambda^\perp$ of a subgroup $\Lambda$? After our discussion of symplectic
characters it should be clear that the correct choice is ${\rm
J}\Lambda^\perp$ because this set consists of all points $(\ulrmk,\ulrml)$
such that ${\rm c}_{|G|}(\lambda,(\ulrmk,\ulrmr))\overline{{\rm
c}_{|G|}}((\ulrmk,\ulrmr),\lambda)=1$ for all $\lambda\in\Lambda$. First this
implies that the set of all points $(\ulrmk,\ulrmr)\in\rmGG$ which satisfies
this condition is another subgroup $\Lambda^\circ$ of $\rmGG$. In Gabor
analysis $\Lambda^\circ$ is called the {\it adjoint subgroup} of $\Lambda$ and
its relevance was first discovered by Feichtinger and Kozek in \cite{feko98}.
Their approach provided an explanation of Janssen's representation of the
frame operator for separable lattices in $\mathbb{R}^{\rm
d}\times\widehat{\mathbb{R}}^{\rm d}$ in \cite{dalala95,ja95} and for certain
non-separable lattices due to Ron and Shen in \cite{rosh93,rosh97}. For
different reasons Rieffel was led to consider $\Lambda^\circ$, which he calls
the orthogonal subgroup of $\Lambda$ \cite{ri88}, and recently Digernes and
Varadarajan came across this object $\Lambda^\circ$ and they refer to it as
the polar of $\Lambda$, \cite{diva04}.

\begin{theorem}[Poisson summation formula]
For any subgroup $\Lambda \lhd \rmGG$  one has:
\begin{equation}\label{Poisson-fin}
  \sum_{\lambda\in\Lambda}{\bf F}(\lambda)=|\Lambda|^{-1}
\sum_{\lambda^\circ\in\Lambda^\circ}{\bf F}(\lambda^\circ) \quad \mbox{for
all} \quad {\bf F} \in \rmCCGG.
 \end{equation}
\end{theorem}
An important consequence of the Poisson summation formula for the symplectic
Fourier transform is the {\it Fundamental Identity of Gabor Analysis}(FIGA).
The following lemma about the symplectic Fourier transform of two STFT's
appears at different places in the engineering and mathematical literature,
see \cite{fasc04} in various degrees of generality. Following \cite{su62} we
call it the Sussman identity. Therefore the next result is a finite analogue
of Sussman's identity.
\begin{proposition}
Let $\Bff_1,\Bff_2,\Bfg_1,\Bfg_2$ be in $\rmCCG$. Then
\begin{equation*}
{\mathcal F}^s\big[{\rm V}_{\Bfg_1}\Bff_1\overline{{\rm
V}_{\Bfg_2}\Bff_2}\big]((\ulrml,\ulrms))={\rm V}_{\Bfg_1}\Bfg_2\overline{{\rm
V}_{\Bff_1}\Bff_2}((\ulrml,\ulrms)).
\end{equation*}
\end{proposition}

\begin{proof}
\begin{eqnarray*}
   {\mathcal F}^s\big[{\rm V}_{\Bfg_1}\Bff_1\overline{{\rm V}_{\Bfg_2}\Bff_2}\big]((\ulrml,\ulrms))&=&\sum_{\ulrmk,\ulrmr}{\rm V}_{\Bfg_1}\Bff_1(\ulrmk,\ulrmr)\overline{{\rm V}_{\Bfg_2}\Bff_2(\ulrmk,\ulrmr)}e^{2\pi i\Omega((\ulrml,\ulrms),(\ulrmk,\ulrmr))} \\
                   &=&\sum_{\ulrmk,\ulrmr}             \langle\pi(\ulrmk,\ulrmr)\Bff_1,\pi(\ulrml,\ulrms)\pi(\ulrmk,\ulrmr)\Bfg_1\rangle\overline{\langle\Bff_2,\pi(\ulrmk,\ulrmr)\Bfg_2\rangle}e^{2\pi i\Omega((\ulrml,\ulrms),(\ulrmk,\ulrmr))}\\
                   &=&\sum_{\ulrmk,\ulrmr}\langle\pi(\ulrml,\ulrms)\Bff_1,\pi(\ulrmk,\ulrmr)\pi(\ulrml,\ulrms)\Bfg_1\rangle\overline{\langle\Bff_2,\pi(\ulrmk,\ulrmr)\Bfg_2\rangle}\\
                   &=&\overline{\langle \Bff_2,\pi(\ulrml,\ulrms)\Bff_1\rangle}\langle\Bfg_2,\pi(\ulrml,\ulrms)\Bfg_1\rangle.
  \end{eqnarray*}
In the last line we have applied Moyal's identity \ref{moyal}.
\end{proof}
The next theorem is just an application of the symplectic Poisson summation
formula to Sussman's Identity, which is the above mentioned {\it Fundamental
Identity of Gabor analysis}.
\begin{theorem}[FIGA]
Let $\Lambda$ be a subgroup of $\rmGG$ and ${\bf f}_1,{\bf f}_2,{\bf
g}_1,{\bf g}_2\in\rmCCG$. Then
\begin{equation*}
  \sum_{\lambda\in\Lambda}{\rm V}_{{\bf g}_1}{\bf f}_1(\lambda)\overline{{\rm V}_{{\bf g}_2}{\bf f}_2}(\lambda)=|\Lambda|^{-1}
\sum_{\lambda^\circ\in\Lambda^\circ}{\rm V}_{{\bf g}_1}{\bf
g}_2(\lambda^\circ)\overline{{\rm V}_{{\bf f}_1}{\bf f}_2}(\lambda^\circ).
 \end{equation*}
\end{theorem}
In the continuous case Janssen was the first to call this identity the FIGA
because important results such as the Wexler-Raz  biorthogonality conditions
and the Ron-Shen duality principle are easily derived from it, \cite{ja95}.
\par
The relevance of the symplectic structure of the time-frequency plane $\rmGG$
was first pointed out by Feichtinger and Kozek in the context of quantization
of operators and Gabor frames for elementary locally compact abelian groups
\cite{feko98}. Their approach invokes the full power of abstract harmonic
analysis and is therefore very technical and not accessible to the majority of workers in Gabor analysis. On the other hand these groups are applying Gabor
analysis to real world problems and one reason for this article is to
transfer the main results of Feichtinger and Kozek into a form which requires
only a modest background in harmonic analysis.

\subsection{Multi-window Gabor frames}

In this section we explore multi-window Gabor frames with the help of the
spreading representation. This allows for a unified discussion of
non-separable Gabor systems which generalize the known results for separable
Gabor systems due to Wexler-Raz, Tolimieri, Qiu and Strohmer,
\cite{rawe90,qi98,st97-3}. First of all we treat just Gabor frames and obtain
their main properties and at the end of this section we indicate how these
results are extended to multi-window Gabor frames, i.e. a finite sum of Gabor
frames.
\begin{definition}\label{def:GabFra}
Let $\Lambda$ be a subgroup of $\rmGG$ and ${\bf g}$ a {\it Gabor atom} in
$\rmCCG$. Then ${\mathcal G}({\bf g},\Lambda)=\{\pi(\lambda){\bf
g}:\lambda\in\Lambda\}$ is called a {\it Gabor system}. If the frame operator
\begin{equation*}
  \rmS_{{\bf g},\Lambda}{\bf f}=\sum_{\lambda\in\Lambda}\langle {\bf f},\pi(\lambda){\bf g}\rangle\pi(\lambda){\bf g}
\end{equation*}
is invertible, then ${\mathcal G}({\bf g},\Lambda)$ is called a {\it Gabor
frame}.
\end{definition}
The signal $\tilde\Bfg:=\rmS_{\Bfg,\Lambda}^{-1}\Bfg$ is the {\it canonical
dual window} and $\Bfh_0:=\rmS_{\Bfg,\Lambda}^{-1/2}\Bfg$ is the {\it
canonical tight window} of the Gabor frame $\mathcal{G}(\Bfg,\Lambda)$.
\par
We consider the elements of the Gabor system
$\{\pi(\lambda)\Bfg:\lambda\in\Lambda\}$ as columns of a
$|G|\times|\Lambda|$-matrix
$\rmDgLa=[\pi(\lambda_1)\Bfg,...,\pi(\lambda_{|\Lambda|})\Bfg]$ for some
ordering of the elements of $\Lambda$. Then the Gabor frame operator
$\rmSgLa$ may be written as $\rmSgLa=\rmDgLa\circ\rmDgLa^\star$, i.e.\
$\rmSgLa$ acts on vectors in $\rmCCG$ and implements the coefficient operator
for ${\mathcal G}(\Bfg,\Lambda)$. Observe that the operator
$\rmDgLa^\star\circ\rmDgLa$ with respect to the canonical basis
$\{\delta_{\lambda,0}:\lambda\in\Lambda\}$ is represented by the {\bf Gram
matrix} $\rmG_{\Bfg,\Lambda}$ of ${\mathcal G}(\Bfg,\Lambda)$ with entries
$\big(\rmG_{\lambda,\mu}=\langle\Bfg_\mu,\Bfg_\lambda\rangle\big)_{\lambda,\mu\in\Lambda}$.
Consequently the Gram matrix $\rmG_{\Bfg,\Lambda}$ acts on vectors in
$\CC^{|\Lambda|}$.
\par
As a motivation of a general Gabor system $\mathcal{G}(\Bfg,\Lambda)$ we
first treat the Gabor frame $\mathcal{G}(\Bfg,{0}\})$, i.e.
$\{\pi(\ulrmk,\ulrmr):(\ulrmk,\ulrmr)\in\rmGG\}$. Then the coefficient
operator and synthesis operator are given by
\begin{equation*}
{\rm C}_{\mathcal
G}{\Bff}=\Big(\langle{\Bff},\ulrmpikr{\Bfg}\rangle\Big)_{(\ulrmk,\ulrmr)}~~\text{and}~~{\rm
D}_{\mathcal G}{\bf
c}=\sum_{(\ulrmk,\ulrmr)\in\rmG\times\widehat\rmG}c(\ulrmk,\ulrmr)\ulrmpikr{\Bfg}~~\text{for}~~{\bf
c}\in\rmCCGG.
\end{equation*}
\begin{proposition}\label{pro:TFtight}
For ${\Bfg}\in\rmCCG$ 
we have that $\mathcal{G}=\{\ulrmpikr{\Bfg}:(\ulrmk,\ulrmr)\in\rmGG\}$ is a
tight frame for $\rmCCG$ with frame constants ${\rm
A,B}=\|\Bfg\|_{\rmCCG}^2$, i.e.
\begin{equation*}
  \|\Bfg\|_{\rmCCG}^2\|{\Bff}\|_{\rmCCG}^2=\frac{1}{|\rmG\times\widehat\rmG|}\sum_{\ulrmk,\ulrmr\in\rmG}|\langle {\Bff},\ulrmpikr\Bfg\rangle|^2~~\text{for all}~~\Bff\in\rmCCG.
\end{equation*}
\end{proposition}
\begin{corollary}[Resolution of Identity]\label{cor:ResId}
Let ${\Bfg}\in\rmCCG$ with $\|{\Bfg}\|_{\rmCCG}=1$. Then for every
${\bf h}\in\rmCCG$ with $ \langle\Bfg,\Bfh\rangle \neq 0 $ one has 
\begin{equation*}
{\bf f}=\frac{1}{|\rmG|\cdot
\langle\Bfg,\Bfh\rangle}\sum_{(\ulrmk,\ulrmr)\in\rmG\times\widehat\rmG}\langle{\Bff},\ulrmpikr\Bfg\rangle\ulrmpikr{\bf
h}.
\end{equation*}
\end{corollary}
The corollary follows from our discussion of frames in the introduction with
${\bf h}_0:=\rmS_{\mathcal G}^{-1}{\Bfg}$. Our proof that
$\{\ulrmpikr{\Bfg}:(\ulrmk,\ulrmr)\in\rmGG\}$ is a tight frame for $\rmCCG$
relies on the commutation relations for time-frequency shifts and Schur's
lemma. The frame operator $\rmS_{\mathcal G}$ of
$\{\ulrmpikr{\Bfg}:(\ulrmk,\ulrmr)\in\rmGG\}$ has the following property:
\begin{lemma}
  For every $(\ulrml,\ulrms)$ in $\rmGG$ we have
\begin{equation*}
  \rmS_{\mathcal G}=\pi(\ulrml,\ulrms)\circ\rmS_{\mathcal G}\circ\pi(\ulrml,\ulrms)^\star.
\end{equation*}
\end{lemma}
\begin{proof}
The commutation relations  $\ulrmpils^\star\ulrmpikr=e^{-2\pi
i(\ulrmk-\ulrml)\cdot\ulrms/|\rmG|}\pi(\ulrmk-\ulrml,\ulrmr-\ulrms)$ yields
\begin{eqnarray*}
   \pi(\ulrml,\ulrms)\circ\rmS_{\mathcal G}\circ\pi(\ulrml,\ulrms)^\star\Bff&=&\sum_{\ulrmk,\ulrmr\in\rmG}\langle\ulrmpils^\star\Bff,\ulrmpikr\Bfg\rangle\ulrmpikr\Bfg\\
  &=&\sum_{\ulrmk,\ulrmr\in\rmG}\langle\Bff,\pi(\ulrmk-\ulrml,\ulrmr-\ulrms)\Bfg\rangle\pi(\ulrmk-\ulrml,\ulrmr-\ulrms)\Bfg=\rmS_{\mathcal G}\Bff.
\end{eqnarray*}
\end{proof}
Recall that $\{\ulrmpikr:(\ulrmk,\ulrmr):\rmGG\}$ is an irreducible
representation of $\rmGG$ and that it generates $\MMG$, hence the preceding
observation yields that $\rmS_{\mathcal G}$ is in the commutant of $\MMG$.
Hence by Schur's lemma $\rmS_{\mathcal G}$ is a multiple of the identity
operator,
\begin{equation*}
  \rmS_{\mathcal G}={\rm C} \cdot \mathbb{I}_{|\rmG|} \quad \text{for some}~~{\rm C}\in\CC.
\end{equation*}
The determination of the constant ${\rm C}$ follows from an application of Moyal's identity. This implies
the assertion that $\{\ulrmpikr:(\ulrmk,\ulrmr)\in\rmGG\}$ is a tight frame
for $\rmCCG$. 
\par
Before we move on to general Gabor frames we note some important properties
of the spreading representation, which allows us to justify the name
"Resolution of Identity". First of all we remark that the projective
representation of $\rmGG$ gives rise to a unitary representation on $\MMG$.
In other words
\begin{equation*}
 ({\ulrmk},{\ulrmr})\mapsto  \rmA[({\ulrmk},{\ulrmr})]=\pi({\ulrmk},{\ulrmr})\rmA\pi({\ulrmk},{\ulrmr})^*,~~\rmA\in\MMG
\end{equation*}
is an involutive automorphism of $\MMG$.
\begin{proposition}
 The mapping  $ ({\ulrmk},{\ulrmr}) \mapsto  \rmA[({\ulrmk},{\ulrmr})]$
 defines a unitary representation of $\rmGG$ on the Hilbert space $\MMG$ with the Frobenius norm,
 i.e.\
 \begin{enumerate}
  \item $\rmA[({\ulrmk},{\ulrmr})]\circ \rmA[({\ulrml},{\ulrms})]=\rmA[{\ulrmk}+{\ulrml},{\ulrmr}+{\ulrms}]$,
  \item $\langle \rmA[({\ulrmk},{\ulrmr})],\rmB[({\ulrmk},{\ulrmr})]\rangle_{\Fro}=\langle \rmA,\rmB\rangle_{\Fro}$.
\end{enumerate}
\end{proposition}
We are interested in the relation between the spreading function of an
operator $\rmA\in\MMG$ and of $\rmA[({\ulrmk},{\ulrmr})]$ for
$({\ulrmk},{\ulrmr})\in\rmGG$. Since the spreading representation of $\rmA$
is an expansion with respect to the basis
$\{\ulrmpikr:({\ulrmk},{\ulrmr})\in\rmGG\}$ the problem reduces to the
understanding of conjugation by $\pi({\ulrml},{\ulrms})$ for time-frequency
shifts.
\begin{lemma}
 For $\rmA\in\MMG$ and $({\ulrmk},{\ulrmr})\in\rmGG$ one has
 $  \eta_{\rmA[({\ulrmk},{\ulrmr})]}={\rm M}^{\bf s}_{({\ulrmk},{\ulrmr})}\eta_{\rmA}$.
\end{lemma}
\begin{proof}
\begin{eqnarray*} \quad \quad \ulrmpikr\circ\rmA\circ\ulrmpikr^\star&=&\sum_{\ulrml,\ulrms}\eta_{\rmA}(\ulrml,\ulrms)\ulrmpikr\ulrmpils\ulrmpikr^\star\\&=&\sum_{\ulrml,\ulrms}\chi_{\rmG}\big(({\ulrmk},{\ulrmr}),({\ulrml},{\ulrms})\big)\overline{\chi}_{\rmG}\big(({\ulrml},{\ulrms}),({\ulrmk},{\ulrmr})\big)\eta_{\rmA}(\ulrml,\ulrms)\ulrmpils\\&=&\sum_{\ulrml,\ulrms}{\rm M}^{\bf s}_{({\ulrml},{\ulrms})}\eta_{\rmA}\ulrmpils.
\end{eqnarray*}
\end{proof}
This behaviour of the spreading coefficients under conjugation of the
operator will be crucial in our discussion of the Janssen representation of
Gabor frame operators. At the moment we want to point out that the preceding
properties of the spreading representation gives the following form of the
Resolution of Identity \eqref{cor:ResId}. Namely, recall the rank-one
operator ${\rm P}_{{\bf g},{\bf h}}{\bf f}=\langle{\bf f},{\bf g}\rangle{\bf
h}$ for $\Bff,\Bfg\in\rmCCG$. Then an elementary calculation gives that
\begin{equation*}
  \big[\ulrmpikr\circ {\rm P}_{{\bf g},{\bf h}}\circ\ulrmpikr^\star\big]\Bff=\langle\Bff,\ulrmpikr\Bfg\rangle\ulrmpikr\Bfh.
\end{equation*}
Therefore conjugation of a rank-one operator ${\rm P}_{{\bf g},{\bf h}}$ by a
time-frequency shift moves it to the point $(\ulrmk,\ulrmr)$ in the
time-frequency plane $\rmCCGG$. We denote conjugation by $\ulrmpikr$ of a
linear operator $\rmA$ on $\rmCCG$ by
\begin{equation*}
  \rmA[(\ulrmk,\ulrmr)]:=\ulrmpikr\circ A\circ\ulrmpikr^\star.
\end{equation*}
Consequently the Resolution of Identity can be expressed in the following
way,
\begin{equation*}
  \mathbb{I}_{|\rmG|}=\frac{1}{\langle\Bfg,\Bfh\rangle}\sum_{\ulrmk,\ulrmr}{\rm P}_{{\bf g},{\bf h}}[(\ulrmk,\ulrmr)],
\end{equation*}
for $\langle\Bff,\Bfh\rangle\ne 0$. If ${\bf g}={\bf h}$ then the identity
operator is a linear combination of orthogonal projections onto the
one-dimensional spaces generated by the family $\{\pi({\ulrmk},{\ulrmr}){\bf
g}\}$. Recall that the pure states of ${\mathcal M}_{|\rmG\times\rmG|}(\CC)$ are the rank one operators. Therefore the resolution of
identity may be understood as shifting a pure state of ${\mathcal
M}_{|\rmG\times \rmG|}(\CC)$ over the discrete time-frequency plane
$\rmCCGG$.
\par
After these preparations  we want  to explore the structure of Gabor frames
${\mathcal G}(\Bfg,\Lambda)$ for
$\Lambda\vartriangleleft\rmG\times\widehat\rmG$, i.e.\ discrete analogues of
the Resolution of Identity:
\begin{equation*}
  \Bff = \sum_{\lambda\in\Lambda} \, \langle\Bff,\pi(\lambda)\Bfh\rangle \,
   \pi(\lambda)\Bfg \quad \text{for a suitable}~~\Bfh\in\rmCCG.
\end{equation*}
The last equation indicates that operators of the type
\begin{equation*}
\rmS_{{\bf g},{\bf h},\Lambda}{\bf f}=\sum_{\lambda\in\Lambda} \langle {\bf
f},\pi(\lambda){\bf h}\rangle\pi(\lambda){\bf g}
\end{equation*}
are closely related with the frame operator of a Gabor system ${\mathcal
G}(\Bfg,\Lambda)$. Due to this fact they are called {\it Gabor frame-type
operators}. The following property of Gabor frame-type operators is crucial
for an understanding of the structure of Gabor systems and is the very reason
for all duality principles in Gabor analysis.
\begin{lemma}
Let $\Lambda$ be a subgroup of $\rmGG$ and ${\bf g}$ a {\it Gabor atom} in
$\rmCCG$. Then the Gabor frame-type operator $\rmS_{{\bf g},{\bf h},\Lambda}$
commutes with $\pi(\lambda)$ for all $\lambda\in\Lambda$, i.e.
\begin{equation}\label{Sgh-commut}
  \pi(\lambda)\circ\rmS_{{\bf g},{\bf h},\Lambda}\circ\pi(\lambda)^\star=\rmS_{{\bf g},{\bf h},\Lambda}.
\end{equation}
\end{lemma}
\begin{proof}
The commutation relation for time-frequency shifts yields for
$\mu\in\Lambda$.:
\begin{eqnarray*}
  [\pi(\mu)\rmS_{{\bf g},{\mathbf{h}},\Lambda}\pi(\mu)^\star] \Bff&=&\sum_{\lambda\in\Lambda}\langle\pi(\mu)^\star{\bf f},\pi(\lambda){\bf h}\rangle\pi(\mu)\pi(\lambda){\bf g}\\
    &=&\sum_{\lambda\in\Lambda}\langle{\bf f},\pi(\mu)\pi(\lambda){\bf h}\rangle\pi(\mu)\pi(\lambda){\bf g}\\
    &=&\sum_{\lambda\in\Lambda}\langle{\bf f},\pi(\mu+\lambda){\bf h}\rangle\pi(\mu+\lambda){\bf g}=\rmS_{{\bf g},{\mathbf{h}},\Lambda}\Bff.
\end{eqnarray*}
\end{proof}
Since  (\ref{Sgh-commut}) is a crucial property of $ \rmS_{{\bf g},{\bf
h},\Lambda}$ 
we want to explore the structure of the set  $\B(\Lambda)$ of all linear
operators $\rmA\in{\mathcal M}_{|\rmG\times\rmG|}$ which commute with all time-frequency
shifts $\pi(\lambda)$ from a given subgroup $\Lambda$ of $\rmGG$. We call
these operators $\Lambda${\it -invariant}. 
\begin{proposition}
For every subgroup  $\Lambda \lhd G\times\widehat{G}$ the set  $\B(\Lambda)$
is an involutive subalgebra of $M_{|\rmG\times\rmG|}(\CC)$.
\end{proposition}
\begin{proof}
  If $\rmA,\rmB$ are $\Lambda$-invariant operators, then we have that:
  \begin{enumerate}
    \item $\rmA+\rmB\in\B(\Lambda)$ and $c\rmA\in\B(\Lambda)$ for all $c\in\CC$.
    \item $\rmA\rmB\in\B(\Lambda)$ since $\pi(\lambda)\rmA\pi(\lambda)^*=\rmA$ and $\pi(\lambda)\rmB\pi(\lambda)^*=\rmB$ implies that
    \begin{equation*}
      \rmA\rmB=\pi(\lambda)\rmA\pi(\lambda)^*\pi(\lambda)\rmB\pi(\lambda)^*=\pi(\lambda)\rmA\rmB\pi(\lambda)^*.
    \end{equation*}
    \item $\rmA\in\B(\Lambda)$ implies $\pi(\lambda)\rmA^\star\pi(\lambda)^*=\rmA^\star$.
  \end{enumerate}
\end{proof}
The algebra $\mathcal{B}(\Lambda)$ of $\Lambda$-invariant operator $\rmA$ is
the commutant of the twisted group algebra $\mathcal{A}(\Lambda)$ within
${\mathcal M}_{|\rmG|}$  is equal to the twisted group algebra
$\mathcal{A}(\Lambda^\circ)$. In detail:
%
\begin{proposition}
Let $\Lambda$ be a subgroup of $\rmGG$. Then we have
\begin{enumerate}
  \item $\mathcal{B}(\Lambda)=\mathcal{A}(\Lambda^\circ)$;
  \item The commutant of $\mathcal{B}(\Lambda)$ is $\mathcal{A}(\Lambda)$ and the commutant of $\mathcal{A}(\Lambda^\circ)$ is $\mathcal{A}(\Lambda)$;
  \item The center of $\mathcal{B}(\Lambda)=\mathcal{A}(\Lambda\cap\Lambda^\circ)$;
  \item $\mathcal{B}(\Lambda)$ is commutative if and only if $\Lambda^\circ\subseteq\Lambda$.
 \end{enumerate}
\end{proposition}
The last assertions indicate that the structure of $\Lambda$ is essential for
the properties of $\mathcal{B}(\Lambda)$. In analogy to symplectic vector
spaces we call the subgroup $\Lambda$ of $\rmGG$ {\it isotropic} if the
symplectic form $\Omega$ vanishes identically on $\Lambda$. The largest
subgroup of $\rmGG$ with this property is naturally called {\it maximal
isotropic}. A moment of reflection shows that a $\Lambda$ is isotropic if and
only if $\Lambda^\circ\subseteq\Lambda$. Therefore $\mathcal{B}(\Lambda)$ is
commutative if and only if $\Lambda$ is isotropic. Now the maximal
commutative subalgebra of $\MMG$ is the algebra of diagonal matrices which
implies that $\Lambda$ is maximal isotropic if and only if $\Lambda$ is a
product lattice $\Lambda\times\Lambda^\perp$ for
$\Lambda\vartriangleleft\rmG$ and
$\Lambda^\perp\vartriangleleft\widehat\rmG$. In other words, the Gabor frame
operator $\rmS_{\mathcal G,\Lambda\times\Lambda^\perp}$ for a product lattice
$\Lambda\times\Lambda^\perp$ is unitarily equivalent to a diagonal operator.
The Zak transform is the interwinding operator which diagonalizes the Gabor
frame operator $\rmS_{\mathcal G,\Lambda\times\Lambda^\perp}$, see
\cite{gr98}.
\par
Note that $\{\pi(\lambda):\lambda\in\Lambda\}$ defines a {\it reducible}
projective representation of $\Lambda$ if  $\Lambda$ is a proper subgroup,
because then its commutant $\mathcal{A}(\Lambda^\circ)$ is non-trivial.
\begin{proposition}
Let $\Lambda$ be a subgroup of $\rmGG$. Then
$\{\pi(\lambda):\lambda\in\Lambda\}$ defines  a reducible projective
representation of $\Lambda$.
\end{proposition}
The previous statement is valid for any subgroup $\Lambda$ of $\rmGG$,
especially for $\Lambda^\circ$. The twisted group algebras
$\mathcal{A}(\Lambda)$ and $\mathcal{A}(\Lambda^\circ)$ are the matrix
algebras underlying Gabor analysis over finite abelian groups. Therefore we
investigate  their structure in  detail. We use the notion of {\it matrix
coefficients for a twisted group algebra}. Let $\pi_{\mathcal{A}(\Lambda)}$
be a representation of the twisted group algebra $\mathcal{A}$. Then we
define a {\it matrix coefficient} of $\mathcal{A}(\Lambda)$ as
\begin{equation*}
  \langle \pi_{\mathcal{A}(\Lambda)}({\bf a})\Bfg,\Bfh\rangle_{\CC^{|\Lambda|}}=\sum_{\lambda\in\Lambda}a(\lambda)\langle\pi(\lambda)\Bfg,\Bfh\rangle_{\CC^{|\Lambda|}}~~\text{for}~~{\bf a}=(a(\lambda)),\Bfg,\Bfh\in\CC^{|\Lambda|}.
\end{equation*}
There is a close relation between the matrix coefficients of
$\mathcal{A}(\Lambda)$ and $\mathcal{A}(\Lambda^\circ)$.
\begin{theorem}
Let $\Lambda$ be a subgroup of $\rmGG$. Then we have
\begin{equation*}
 \sum_{\lambda\in\Lambda}\sum_{\rm i=1}^{\rm n}{\rm d_i}\langle\Bfg_{\rm i},\pi(\lambda)\Bfh_{\rm i}\rangle\langle\pi(\lambda)\Bfg,\Bfh\rangle=\frac{1}{|\Lambda|}\sum_{\lambda^\circ\in\Lambda^\circ}\sum_{\rm i=1}^{\rm n}{\rm d_i}\langle\Bfg_{\rm i},\pi(\lambda^\circ)\Bfg\rangle\langle\pi(\lambda^\circ)\Bfh,\Bfh_{\rm i}\rangle.
\end{equation*}
\end{theorem}
\begin{proof}
Recall that every $\rmA\in\mathcal{A}(\Lambda)$ may be written as
\begin{equation*}
\pi_{\mathcal{A}(\Lambda)}({\eta_{\rmA}})=\sum_{\lambda\in\Lambda}\sum_{\rm
i=1}^{\rm n}{\rm d_i}\langle\Bfg_{\rm i},\pi(\lambda)\Bfh_{\rm
i}\rangle\pi(\lambda)
\end{equation*}
 and this yields to
\begin{eqnarray*}
 \langle\pi_{\mathcal{A}(\Lambda)}({\eta_{\rmA}})\Bfg,\Bfh\rangle_{\CC^{|\Lambda|}}&=&\sum_{\lambda\in\Lambda}\sum_{\rm i=1}^{\rm n}{\rm d_i}\langle\Bfg_{\rm i},\pi(\lambda)\Bfh_{\rm i}\rangle\langle\pi(\lambda)\Bfg,\Bfh\rangle\\
   &=&\frac{1}{|\Lambda|}\sum_{\lambda^\circ\in\Lambda^\circ}\sum_{\rm i=1}^{\rm n}{\rm d_i}\langle\Bfg_{\rm i},\pi(\lambda^\circ)\Bfg\rangle\langle\pi(\lambda^\circ)\Bfh,\Bfh_{\rm i}\rangle,
\end{eqnarray*}
where we applied in the last equation the FIGA. The right side of the last
equation may be understood as the matrix coefficient for a certain element of
$\pi_{\mathcal{A}(\Lambda^\circ)}$.
\end{proof}
If we consider the special case $\Lambda=\rmGG$, then the statement of the
last theorem specializes to Moyal's Formula, i.e. the Schur orthogonality
relations for STFT's. In this sense we consider the last theorem as a
generalization of Schur's orthogonality relation to reducible group
representations.
\par
The preceding observations yield the following representation of a
$\Lambda$-invariant operator, which includes the representation of Gabor
frame operators as given by Tolimieri-Orr, Qiu, Strohmer, and Wexler-Raz for
product lattices $\Lambda_1\times\Lambda_2$.
\begin{theorem}[Janssen representation]
Let $\Lambda$ be a subgroup of $\rmGG$. Then for a $\Lambda$-invariant
operator $\rmA$ we have a {\it prototype matrix} ${\rm P}$ such that
\begin{equation*}
  \rmA=\sum_{\lambda\in\Lambda}{\rm P}[\lambda]=\sum_{\lambda\in\Lambda} \,  \pi(\lambda)^\star\circ{\rm
  P} \circ \pi(\lambda)
\end{equation*}
is the $\Lambda$-periodization of ${\rm P}$ or
\begin{equation*}
  \rmA=\sum_{\lambda^\circ\in\Lambda^\circ}\langle{\rm P},\pi(\lambda^\circ)\rangle_{\operatorname{Fro}}\pi(\lambda^\circ),
\end{equation*}
i.e. the spreading coefficients of $\rmA$ are the sampled spreading
coefficients of the prototype operator ${\rm P}$ to the adjoint lattice
$\Lambda^\circ$.
\end{theorem}
\begin{proof}
An application of the commutation relations for time-frequency shifts in the
following form
\begin{equation*}
  \pi(\mu)\pi(\lambda)\pi(\mu)^*={\rm c}_{|\rmGG|}^{\bf s}(\lambda,\mu)\pi(\lambda-\mu) \quad \text{for all}~~\lambda,\mu\in\Lambda
\end{equation*}
implies that $\pi(\mu)\circ\pi_{\mathcal{A}(\Lambda)}({\bf a})\circ\pi(\mu)$
is another element of $\mathcal{A}(\Lambda)$ for a suitable translated and
shifted version of ${\bf a}$. In terms of the spreading representation, the imposition of $\Lambda$-invariance on a linear operator $\rmA$ implies that the spreading coefficients are periodic,i.e. $\rmA\in\B(\Lambda)$ is
equivalent to $\rmA[\lambda]=\rmA$, i.e.\ that for all $\lambda \in\Lambda$
one has:
\begin{equation*}
  \eta_{\rmA[\lambda]}(\mu)={\rm c}_{|G|}(\lambda,\mu)\overline{{\rm c}_{|G|}}
  (\mu,\lambda)\eta_{\rmA}(\mu)~~\text{for all}~~\mu\in\Lambda.
\end{equation*}
In other words every $\Lambda$-invariant operator has a representation of the
form
\begin{equation*}
  \rmA=\sum_{\lambda^\circ\in\Lambda^\circ}{\rm a}(\lambda^\circ)\pi(\lambda^\circ)
\end{equation*}
for some vector ${\bf a}=\big(a(\lambda^\circ)\big)\in\CC^{|\Lambda^\circ|}$.
\par
The prototype operator ${\rm P}$ is the sum of projection operators arising
from the singular value decomposition of $\rmA$, i.e.\  ${\rm P}=\sum_{\rm
i=1}^{\rmr}{\rm d_i}{\bf u}_{\rm i}\otimes\ol{{\bf v}_{\rm i}}$. Let $\rmA$
be the Gabor frame operator of ${\mathcal G}(\Bfg,\Lambda)$. Then
\begin{equation*}
  \rmS_{\Bfg,\Lambda}=|\Lambda|^{-1}\sum_{\lambda^\circ\in\Lambda^\circ}\langle \Bfg,\pi(\lambda^\circ)\Bfg\rangle\pi(\lambda^\circ).
\end{equation*}
\end{proof}
We interpret the Janssen representation of a $\Lambda$-invariant operator in
terms of multi-window Gabor frames. At the moment we draw some consequences
of the Janssen representation which generalize results of Wexler-Raz on Gabor
frame operators for separable lattices in cyclic groups. In their work on the
inversion of Gabor frame operators over cyclic groups Wexler and Raz gave the
impetus for the duality theory of Gabor frames which was developed
independently by several groups of researchers. Our general results on the
structure of $\Lambda$-invariant operators lead naturally to the study of a
finite number of Gabor frames. Let ${\bf g}_1,...,{\bf g}_{\rm r}\in\rmCCG$ be Gabor
atoms of the Gabor systems ${\mathcal G}({\bf g}_1,\Lambda),...,{\mathcal
G}({\bf g}_{\rm r},\Lambda)$, then we speak of a {\it multi-window Gabor
system}. If the associated frame operator
\begin{equation*}
  \rmS_{\mathcal G}{\bf f}=\sum_{\rm j=1}^{\rm r}\rmS_{{\bf g}_{\rm j},\Lambda}f=\sum_{\rm j=1}^{\rm r}\sum_{\lambda\in\Lambda}\langle {\bf f},\pi(\lambda){\bf g}_{\rm j}\rangle\pi(\lambda){\bf g}_{\rm j}
\end{equation*}
is invertible, then the system ${\mathcal G}=\bigcup_{j=1}^r{\mathcal G}({\bf
g}_j,\Lambda)$ is called a {\it multi-window Gabor frame}. Therefore our
discussion of the Janssen representation of a $\Lambda$-invariant operator
$\rmA$ states that $\rmA$ is a multi-window Gabor frame operator with
$\Bfg_{\rm i}={\bf u}_{\rm i}$ for ${\rm i}=1,...,\rmr$. We refer the reader
to \cite{zezi93,pozezi98,dofegr06} for further information on multi-window
Gabor frames.

The multi-window Gabor frame operator is actually a finite sum of rank-one
operators
\begin{equation*}
  \rmS_{\mathcal G}=\sum_{\lambda\in\Lambda}({\bf g}_{\rm 1}\otimes\ol{{\bf g}_{\rm 1}}[\lambda]+\cdots+{\bf g}_{\rm r}\otimes\ol{{\bf g}_{\rm r}}[\lambda])
\end{equation*}
and is a $\Lambda$-invariant operator and therefore has a Janssen
representation.

\par
The main problem in Gabor analysis over finite groups is an understanding of
all dual pairs $(\Bfg,\Bfh)$ for a given Gabor system ${\mathcal
G}(\Bfg,\Lambda)$. In other words, we look for all $(\Bfg,\Bfh)$ such that
$\rmS_{\Bfg,\Bfh,\Lambda}=\mathbb{I}_{\rmG}$. We attack this problem by an
application of the Janssen representation of $\Lambda$-invariant operators.
\begin{corollary}
Let $\Lambda$ be a subgroup of $\rmGG$ and ${\mathcal G}(\Bfg,\Lambda)$ a
Gabor system for $\Bfg\in\rmCCG$. Then the following are equivalent:
\begin{enumerate}
  \item $\Bff=\sum_{\lambda\in\Lambda}\langle\Bff,\pi(\lambda)\Bfh\rangle\pi(\lambda)\Bfg \, \mbox{ for all} \, \,  \bf f\in\rmCCG$ .
  \item $\langle\Bfg,\pi(\lambda^\circ)\Bfh\rangle=|\Lambda| \cdot \delta_{\lambda^\circ,\, 0}$ for all $\lambda^\circ\in\Lambda^\circ$.
\end{enumerate}
\end{corollary}
\begin{proof}
The equivalence of (1) and (2) follows immediately from the uniqueness of the
spreading representation, the fact that  $\eta_{\,
\mathbb{I}_{|\rmG|}}=\delta_{\lambda^\circ,\, 0}$, and the identity
  \begin{equation*}
    S_{\Bfg,\Bfh,\Lambda}\Bff=|\Lambda|^{-1}\sum_{\lambda^\circ\in\Lambda^\circ}\langle \Bfg,\pi(\lambda^\circ)\Bfh\rangle\pi(\lambda^\circ)\Bff,
\, \mbox{ for all} \, \,  \bf f\in\rmCCG.
  \end{equation*}
%
\end{proof}
The biorthogonality relations of Wexler-Raz are a relation between traces of
$\mathcal{A}(\Lambda)$ and $\mathcal{A}(\Lambda^\circ)$. Namely, the
commutation relations for time-frequency shifts and the tracial property
yields that
\begin{equation*} \operatorname{tr}_{\mathcal{A}(\Lambda)}(\rmA)=\eta_{\rmA}(0)=a_{0,0}~~\text{and}~~\operatorname{tr}_{\mathcal{A}(\Lambda^\circ)}(\rmA)={\rm c}\eta_{\rmA}(0)~~\text{for some}~~~~{\rm c}\in\CC.
\end{equation*}
The Poisson summation formula for the symplectic Fourier transform allows us to
identity ${\rm c}$ with $1/|\Lambda|$, i.e. the traces of
$\mathcal{A}(\Lambda)$ and $\mathcal{A}(\Lambda^\circ)$ are multiples of each
other:
\begin{equation*}
 \operatorname{tr}_{\mathcal{A}(\Lambda)}=\frac{1}{|\Lambda|}\operatorname{tr}_{\mathcal{A}(\Lambda^\circ)}.
\end{equation*}
Let $\rmA$ be a Gabor frame-type operator. Then the preceding equation yields
the Wexler-Raz biorthogonality relations. More generally, we get a "weighted"
Wexler-Raz biorthogonality relation for multi-window Gabor frame-type
operators.
\par
The Janssen representation of a Gabor frame operator provides an alternative
route to the inversion of the Gabor frame operator: Namely, a Gabor
frame operator may be considered as a twisted convolution operator on
$\Lambda^\circ$. More generally,
\begin{lemma}
  Let $\Lambda$ be a subgroup of $\rmGG$ and $\rmA,\rmB\in\B(\Lambda)$. Then the spreading coefficients of $\rmA\rmB$ are given by the twisted convolution of those for $\rmA$ and $\rmB$:
  \begin{equation*} \eta_{\rmA\rmB}(\lambda^\circ)=\sum_{\mu^\circ\in\Lambda^\circ}\eta_{\rmA}(\lambda^\circ-\mu^\circ)\eta_{\rmB}(\mu^\circ){\rm c}_{\rmG}(\lambda^\circ-\mu^\circ,\mu^\circ).
  \end{equation*}
\end{lemma}
Therefore we have a symbolic calculus for $\Lambda$-invariant operators. We
have that
$\rmS_{\Bfg,\Lambda}\times\mathbb{I}_{\Lambda}=\rmS_{\Bfg,\Lambda}$, i.e.\
\begin{equation*} \eta_{\rmS_{\Bfg,\Lambda}}(\lambda^\circ)=\eta_{\rmS_{\Bfg,\Lambda}}\natural_{\Lambda^\circ}\eta_{\mathbb{I}_{\Lambda}}(\lambda^\circ)=|\Lambda|^{-1}\sum_{\mu^\circ}\langle \Bfg,\pi(\lambda^\circ-\mu^\circ)\Bfg\rangle{\rm c}_{\rmG}(\lambda^\circ-\mu^\circ,\mu^\circ).
\end{equation*}
On the level of representation coefficients the discussion of Gabor frames
corresponds to the study of the mapping
\begin{equation*}
  {\bf c}=({\rm c}_{\lambda^\circ})\mapsto {\rmG}{\bf c}
\end{equation*}
with
$\rmG=\big(\rmG(\lambda^\circ,\mu^\circ\in\Lambda^\circ)\big)=\big(\langle\Bfg,\pi(\lambda^\circ-\mu^\circ)\Bfg\rangle{\rm
c}_{\rmG}(\lambda^\circ-\mu^\circ,\mu^\circ)\big)$. By the definition of the
Gram matrix and the commutation relation, $\rmG$ is the Gram matrix of
$\{\pi(\lambda^\circ)\Bfg:\lambda^\circ\in\Lambda^\circ\}$. Now the question
of invertibility of the Gabor frame operator translates into the
invertibility of the twisted convolution of the Gram matrix of ${\mathcal
G}(\Bfg,\Lambda^\circ)$. More precisely, ${\mathcal G}(\Bfg,\Lambda)$ is a
Gabor frame for $\rmCCG$ if and only if $\rmG$ is invertible. Let ${\bf c}$
be a vector in $\CC^{|\Lambda^\circ|}$. Then ${\bf c}$ represents the
canonical dual window $\Bfh_0={\rm S}_{\Bfg,\Lambda}^{-1}\Bfg$, if it is the
unique solution of the system of equations
\begin{equation*}
 \rmG{\bf c}=|\Lambda|\delta_{\lambda^\circ,0}.
\end{equation*}
Consider the algebra of all $\Lambda$-invariant operators $\B(\Lambda)$ as a
$\CS$-algebra of matrices, i.e.\  we equip ${\mathcal
M}_{|\rmG|\times|\Lambda|}$ with the operator norm. Then we may express an
element of $\B(\Lambda)$ in terms of $\{\pi(\lambda):\lambda\in\Lambda\}$ or
$\{\pi(\lambda^\circ):\lambda^\circ\in\Lambda^\circ\}$. Now by the spreading
representation we obtain a relation between the operator norm of
$\|\rmS_{\Bfg,\Lambda}\|_{\operatorname{op}}$ and
$\|\rmG\|_{\operatorname{op}}$. By the uniqueness of the $\CS$-algebra norm
we get the existence of a constant $\rm k$ such that
$\|\rmS_{\Bfg,\Lambda}\|_{\operatorname{op}}={\rm
k}\|\rmG\|_{\operatorname{op}}$ which for
$\|\rmS_{\Bfg,\Lambda}\|_{\operatorname{op}}$ equal to the identity yields
${\rm k}=|\Lambda|^{-1}$. As a summary of the previous statements we get the
other duality result for Gabor frames, the {\it Ron-Shen duality principle}.
\begin{theorem}[Ron-Shen duality]
Let $\Lambda$ be a subgroup of $G\times\widehat G$. Then ${\mathcal
G}(\Bfg,\Lambda)$ is a frame for $\rmCCG$ if and only if ${\mathcal
G}(\Bfg,\Lambda^\circ)$ is a (Riesz) basis for $\CC^{|\Lambda^\circ|}$. Let
$A_\Lambda, B_\Lambda$ be the frame bounds of ${\mathcal G}(\Bfg,\Lambda)$
and $A_{\Lambda^\circ},B_{\Lambda^\circ}$ the (Riesz) basis bounds of
${\mathcal G}(\Bfg,\Lambda^\circ)$. Then
\begin{equation*}
  A_{\Lambda^\circ}=|\Lambda|A_{\Lambda}~~~\text{and}~~B_{\Lambda^\circ}=|\Lambda|B_{\Lambda}.
\end{equation*}
\end{theorem}
\section{Dual windows,  canonical tight windows and L\"owdin orthogonalization}
In this section we want to stress the importance of the canonical dual and
canonical tight Gabor window, especially the relation between the structure
of the canonical tight Gabor frame and the L\"owdin orthogonalization of the
original Gabor Riesz basis. 
\par
Let ${\mathcal G}(\Bfg,\Lambda)$ be a Gabor frame for $\rmCCG$. Then we
denote by $\Gamma_{\Bfg,\Lambda}$ the set of all dual pairs $(\Bfg,\Bfh)$,
i.e.
\begin{equation*}
  \Gamma_{\Bfg,\Lambda}=\{\Bfh:\rmC_{\Bfg,\Lambda}^\star\rmC_{\Bfg,\Lambda}=\mathbb{I}_{|\rmG|}\}.
\end{equation*}
The exploration of this set  is one of the core problems in Gabor analysis.
Observe that the $\Lambda$-invariance has as important consequence, the
$\Lambda$-invariance of arbitrary powers of the Gabor frame operator:
\begin{equation}\label{E:PowersFrameOperator}
    \pi(\lambda)\circ\rmS_{\Bfg,\Lambda}^\alpha \circ\pi(\lambda)^\star=\rmS_{\Bfg,\Lambda}^\alpha~~\text{for all}~~\alpha\in\mathbb{R}.
\end{equation}
Therefore
\begin{equation*}
\Bff=\rmS_{\Bfg,\Lambda}^{-1}\rmS_{\Bfg,\Lambda}\Bff=\sum_{\lambda\in\Lambda}\langle\Bff,\pi(\lambda)\Bfg\rangle\pi(\lambda)(\rmS_{\Bfg,\Lambda}^{-1}\Bfg)
\end{equation*}
and $\mathcal{G}(\tilde\Bfg,\Lambda)$ is called the {\it canonical dual Gabor
frame} given by $\tilde\Bfg=\rmS^{-1}_{\Bfg,\Lambda} \Bfg$. The duality
theory of Gabor frames implies that
\begin{equation*}
  \rmS_{\tilde g,\Lambda}=\frac{1}{|\Lambda|}\sum_{\lambda^\circ\in\Lambda^\circ}\langle\tilde g,\pi(\lambda^\circ)\tilde g\rangle\pi(\lambda^\circ),
\end{equation*}
i.e. $\rmS_{\tilde g,\Lambda}$ is an element of the linear span of
$\mathcal{G}(\Bfg,\Lambda^\circ)=\{\pi(\lambda^\circ)\Bfg:\lambda^\circ\in\Lambda^\circ\}$.
By the Wexler-Raz biorthogonality relation we can identify $\Gamma_{\Bfg,\Lambda}$.
\begin{theorem}[Wexler-Raz]
Let ${\mathcal G}(\Bfg,\Lambda)$ be a Gabor frame for $\rmCCG$. Then $\Bfh\in
\Gamma_{\Bfg,\Lambda}$ if and only if
\begin{equation*}
  \Bfh\in\tilde\Bfg+\mathcal{G}(\Bfg,\Lambda^\circ)^\perp.
\end{equation*}
\end{theorem}
\begin{proof}
Our argument closely follows \cite{gr01}. Every
$\Bfh\in\Gamma_{\Bfg,\Lambda}$ satisfies the Wexler-Raz biorthogonality
conditions, therefore we have
$\langle\Bfh-\tilde\Bfg,\pi(\lambda^\circ)\Bfg\rangle=0$ for all
$\lambda^\circ\in\Lambda^\circ$, i.e.
$\Bfh-\tilde\Bfg\in\mathcal{G}(\Bfg,\Lambda^\circ)^\perp$. Conversely, letting
$\Bfh\in\tilde\Bfg+\mathcal{G}(\Bfg,\Lambda^\circ)^\perp$, then $\Bfh$
fulfills the Wexler-Raz biorthogonality relation, i.e.
$\rmS_{\Bfg,\Bfh,\Lambda}=\mathbb{I}_{|\rmG|}$.
\end{proof}
The canonical dual Gabor atom $\tilde\Bfg=\rmS^{-1}_{\Bfg,\Lambda}\Bfg$ is a
distinguished element of $\Gamma_{\Bfg,\Lambda}$ since it is the "optimal"
dual window in various ways.
\begin{theorem}
Let ${\mathcal G}(\Bfg,\Lambda)$ be a Gabor frame for $\rmCCG$. Then for all
$\Bfh\in\Gamma_{\Bfg,\Lambda}$ the following (equivalent, characteristic
properties) hold true:
\begin{eqnarray*}
   \|\tilde\Bfg\|_{\rmCCG}&\le&\|\Bfh\|_{\rmCCG} \quad (\text{minimal norm}),\\
   \|\rmC_{\tilde\Bfg,\Lambda}\Bff\|_{\rmCCG}&\le&\|\rmC_{\tilde\Bfh,\Lambda}\Bff\|_{\rmCCG} \quad (\text{minimal norm coefficients}),\\
  \|\Bfg-\tilde\Bfg\|_{\rmCCG}&\le&\|\Bfg-\Bfh\|_{\rmCCG}\quad(\text{closest to Gabor atom}),\\
  \|\Bfg/\|\Bfg\|_{\rmCCG}-\tilde\Bfg/\|\tilde\Bfg\|_{\rmCCG}\|_{\rmCCG}&\le&\|\Bfg/\|\Bfg\|_{\rmCCG}-\Bfh/\|\Bfh\|_{\rmCCG}\|_{\rmCCG} \quad (\text{most likely}).
\end{eqnarray*}
\end{theorem}
We refer the interested reader to \cite{pr96} for the proof of the preceding
theorem and for a detailed discussion of the structure of
$\Gamma_{\Bfg,\Lambda}$ and that it actually characterizes the dual Gabor
atom. We only want to focus on the minimal norm property of the coefficient
mapping for the canonical dual Gabor frame. The following computation  shows
the relevance of the Moore-Penrose inverse for Gabor frames.
\begin{eqnarray*} \mathcal{C}_{\tilde\Bfg,\Lambda}\Bff&=&(\langle\Bff,\pi(\lambda)\rmS_{\Bfg,\Lambda}^{-1}\Bfg\rangle)=(\langle\Bff,\rmS_{\Bfg,\Lambda}^{-1}\pi(\lambda)\Bfg\rangle)\\
&=&(\langle\rmS_{\Bfg,\Lambda}^{-1}\Bff,\pi(\lambda)\Bfg\rangle)=\mathcal{C}_{\tilde\Bfg,\Lambda}(\mathcal{C}_{\tilde\Bfg,\Lambda}^\star\mathcal{C}_{\tilde\Bfg,\Lambda})^{-1}\Bff=\mathcal{C}_{\Bfg,\Lambda}^+\Bff,
\end{eqnarray*}
where $\mathcal{C}_{\Bfg,\Lambda}^+$ denotes the Moore-Penrose inverse of the
coefficient operator $\mathcal{C}_{\Bfg,\Lambda}$. Then it is well-known that
$\mathcal{C}_{\Bfg,\Lambda}^+$ is the minimal solution of the system
$\mathcal{C}_{\Bfg,\Lambda}\Bff={\bf a}$ for a given coefficient vector ${\bf
a}$ and such that $\Bff$ may be expressed in terms of
$(\langle\Bff,\pi(\lambda)\Bfg\rangle)$, i.e. $\Bff={\rm B}{\bf a}$ under the
constraint that
\begin{equation*}
  \|\mathcal{C}_{\Bfg,\Lambda}\Bff-{\bf a}\|_{\CC^{|\rm G|}}^2=\|(\mathcal{C}_{\Bfg,\Lambda}{\rm B}-\mathbb{I}_{|\rmG|}){\bf a}\|_{\CC^{|\rm G|}}^2=\text{minimal}.
\end{equation*}
In other words we have to minimize $\|\mathcal{C}_{\Bfg,\Lambda}{\rm
B}-\mathbb{I}_{|\rmG|}\|_{\operatorname{Fro}}^2$ over all ${\rm B}$ which is
solved by the Moore-Penrose inverse $B={\rm C}_{\Bfg,\Lambda}^+$. In
\cite{gogo00} it has been shown that the Moore-Penrose inverse minimizes
\begin{equation*}
  \|{\rm C}_{\Bfg,\Lambda}{\rm B}-I\|_{\rm S_p}
\end{equation*}
for all Schatten-von Neumann norms, $ \|.\|_{\rm S_p}$ for all
$p\in[1,\infty]$. Consequently, the coefficients of the canonical dual Gabor
frame are minimal for all $\ell^{\rm p}$-norms.
\par
There exists a canonical way to associate a tight Gabor frame with
$\mathcal{G}(\Bfg,\Lambda)$. By the $\Lambda$-invariance of
$\rmS_{\Bfg,\Lambda}$ we have
\begin{equation*} \rmS^{-1/2}_{\Bfg,\Lambda}\rmS_{\Bfg,\Lambda}\rmS^{1/2}_{\Bfg,\Lambda}=\sum_{\lambda\in\Lambda}\langle\Bff,\pi(\lambda)\rmS^{-1/2}_{\Bfg,\Lambda}\Bfg\rangle\pi(\lambda)\rmS^{-1/2}_{\Bfg,\Lambda}\Bfg
\end{equation*}
Observe that ${\mathcal G}(\Bfh^\circ,\Lambda)$ is a tight frame generated by
the canonical tight Gabor atom $\Bfh^\circ=\rmS^{-1/2}_{\Bfg,\Lambda}\Bfg$.
More generally, we investigate the set of all tight Gabor frames ${\mathcal
G}(\Bfh,\Lambda)$ for a given subgroup $\Lambda$ in $\rmCCGG$, i.e.
$\Gamma_{\operatorname{tight}}(\Bfg,\Lambda)=\{\Bfh\in\CC^{|\rmG|}:
\rmS_{{\bf h},\Lambda}={\mathbb{I}_{|{\rm G}|}}\}$.
\begin{theorem}
 Assume that  $\Bfh\in\Gamma_{\text{tight}}(\Bfg,\Lambda)$. Then there exists a unitary operator ${\rm U}$ such that ${\rm U}{\rm C_{\bf h,\Lambda}}={\rm C_{\bf h^\circ,\Lambda}}$.
\end{theorem}
\begin{proof}
  The statement of the theorem may be expressed as follows. If ${\rm C}_{\bf h^\circ,\Lambda}^\star{\rm C}_{\bf h^\circ,\Lambda}={\rm C}_{\bf h,\Lambda}^\star{\rm C}_{\bf h,\Lambda}$ holds,  then there exists a unitary ${\rm U}$ such that ${\rm C}_{\bf h^\circ,\Lambda}={\rm U}{\rm C}_{\bf h^\circ,\Lambda}$. We define ${\rm U}$ on the range ${\text{ran}}({\rm C}_{\bf h^\circ,\Lambda})$ by ${\rm U}{\rm C}_{\bf h^\circ,\Lambda}\Bff={\rm C}_{\bf h,\Lambda}\Bff$ and by zero on the orthogonal complement. Then ${\rm C}_{\bf h^\circ,\Lambda}^\star{\rm C}_{\bf h^\circ,\Lambda}={\rm C}_{\bf h,\Lambda}^\star{\rm C}_{\bf h,\Lambda}$ implies $\text{ker}({\rm C}_{\bf h^\circ,\Lambda})\subseteq\text{ker}({\rm C}_{\bf h,\Lambda})$. Consequently ${\rm U}$ is an isometry satisfying
  ${\rm U}{\rm C}_{\bf h,\Lambda}={\rm C}_{\bf h^\circ,\Lambda}$.
\end{proof}
The previous theorem corresponds to the fact that there is some freedom in
taking square-roots of a positive operator. The canonical tight window
$\Bfh^\circ$ is the tight frame which minimizes $\|\Bfg-\Bfh\|_{\CC^{|G|}}$
among all tight frames. By the Ron-Shen duality principle this is equivalent
to the minimization over all orthonormal basis ${\mathcal G}(\Bfh,\Lambda)$ for
the linear span of ${\mathcal G}(\Bfg,\Lambda)$. The problem for the
orthonormal basis was first solved by L\"owdin in \cite{lo50} and corresponds
to the choice ${\rm A}={\rmG}^{-1/2}=({\rm C}_{\bf g,\Lambda}^\star{\rm
C}_{\bf g,\Lambda})^{-1/2}$. We formulate L\"owdin's orthogonalization for a
general set of linearly independent vectors
$\mathcal{G}=\{\Bfg_1,...,\Bfg_{\rm N}\}$ of $\CC^{\rm M}$. Let ${\rm C}$ be
the ${\rm M}\times{\rm N}$ matrix with $\{\Bfg_1,...,\Bfg_{\rm N}\}$ as
columns.
\begin{theorem}[L\"owdin]
Let $\mathcal{G}=\{\Bfg_1,...,\Bfg_{\rm N}\}$ be a collection of linearly
independent vectors in $\CC^{\rm M}$ and singular value decomposition ${\rm
C}={\rm U}{\rm D}{\rm V}^\star$ of ${\rm C}$. Then $\|{\rm
C}-X\|_{\operatorname{Fro}}$ is minimized for unitary $X$ by ${\rm L}={\rm
C}({\rm C}^\star{\rm C})^{-1/2}={\rm UV}^\star$.
\end{theorem}
\begin{proof}
  First we show that ${\rm L}$ is unitary and secondly we show that it is the unique minimizer. Since ${\rm C}={\rm U}{\rm D}{\rm V}^\star$ we have that
  \begin{eqnarray*}
    {\rm L}&=&{\rm U}{\rm D}{\rm V}^\star({\rm V}{\rm D}{\rm U}^\star{\rm U}{\rm D}{\rm V}^\star)^{-1/2}={\rm U}{\rm D}{\rm V}^\star({\rm V}{\rm D}^2{\rm V}^\star)^{-1/2}\\
        &=&{\rm U}{\rm D}{\rm V}^\star({\rm V}{\rm D}{\rm V}^\star)^{-1}={\rm UV}^\star.
  \end{eqnarray*}
  Since ${\rm U,V}$ are unitary, we have that ${\rm L}$ is unitary.
  \begin{equation*}
    \|{\rm C}-{\rm X}\|_{\operatorname{Fro}}=\|{\rm U}{\rm D}{\rm V}^\star-{\rm X}\|_{\operatorname{Fro}}=\|{\rm D}-{\rm U}^\star{\rm X}{\rm V}\|_{\operatorname{Fro}}=\text{min},
  \end{equation*}
  if ${\rm U}^\star{\rm X}{\rm V}=\mathbb{I}$, i.e. ${\rm X}={\rm U}{\rm V}^\star$.
\end{proof}
In \cite{gole91} it is demonstrated that ${\rm L}$ is the minimal
orthogonalization for any unitarily invariant norm, especially for all
Schatten-von Neumann classes $\|.\|_{\rmS^p}$ and $p\in[0,\infty]$. Many
authors call the L\"owdin orthogonalization the {\it symmetric
orthogonalization} because it is invariant under a permutation of the
linearly independent vectors of ${\rm G}$. More precisely, let ${\rm P}$ be a
${\rm N}\times{\rm M}$ permutation matrix. Then ${\rm LP}$ is the best
orthogonalization of ${\rm C P}$. The proof is by contradiction.
\par
 We close with a few words on the case of Gabor frames. Since the Gram matrix of a Gabor frame
has a certain structure, the absolute value of entries are circulant, and the
Frobenius norm becomes the norm of the Gabor atom, see Janssen and Strohmer
for a more detailed discussion \cite{jast02}.

\section{Conclusion}\label{Conclusion}
The main goal of this paper is both a self-contained description and a survey  of Gabor frames for finite-dimensional
Hilbert spaces including its duality theory, such as the Wexler-Raz biorthogonality
relations, Janssen's representation of the Gabor frame operator and the
Ron-Shen's duality principle. The presentation relies on elementary
properties of twisted group algebras, which are matrix algebras in the given
 finite-dimensional setting. Our focus on twisted group algebras allows
an elementary and self-contained approach. The interplay of these twisted
group algebras underlies the well-known results derived by Wexler-Raz, Qiu,
Strohmer and others on Gabor frames for product lattices. Therefore we are
able to generalize their results to arbitrary subgroups of the time-frequency
plane using algebraic methods instead of complicated multi-index
calculations. The understanding of these twisted group algebras is also
intimately related to the symplectic structure of the time-frequency plane.
Consequently we have devoted a part of this survey to describe harmonic
analysis on the time-frequency plane. As our note is mostly addressed to
applied mathematicians as well as engineers we did not go so far as to
explain in which sense the Morita equivalence of the twisted group algebras
of the subgroup and its adjoint subgroup. However, from a deeper
view-point it can be seen as  the very reason for the validity of duality
theory of Gabor frames (see \cite{lu05-1} for details).


\flushleft{

}

\end{document}